\documentclass[a4paper,12pt]{article}
\usepackage{graphicx}
\usepackage{psfrag}

\usepackage{mathrsfs}
\usepackage{amssymb}
\usepackage{amsfonts}
\usepackage{latexsym}
\usepackage{amsmath}
\usepackage{amsthm}
\usepackage[cp1251]{inputenc}
\usepackage[english]{babel}

\newcommand{\er}[1]{\textrm{(\ref{#1})}}
\def\lb{\label}

\theoremstyle{plain}
\newtheorem{Def}{\bf Definition}
\newtheorem{theorem}{\bf Theorem}[section]
\newtheorem{lemma}[theorem]{\bf Lemma}
\newtheorem{proposition}[theorem]{\bf Proposition}

\theoremstyle{remark}

\renewcommand{\a}{\alpha}           

\renewcommand{\b}{\beta}            

\newcommand{\g}{\gamma}           \newcommand{\cC}{\mathcal{C}}

\newcommand{\G}{\Gamma}           

\renewcommand{\d}{\delta}           

\newcommand{\D}{\Delta}           

\newcommand{\ve}{\varepsilon}     

\newcommand{\z}{\zeta}            \newcommand{\cH}{\mathcal{H}}

\newcommand{\e}{\eta}             

\newcommand{\vt}{\vartheta}       

\newcommand{\vT}{\Theta}          

\renewcommand{\k}{\kappa}           

\renewcommand{\l}{\lambda}          \newcommand{\cM}{\mathcal{M}}

\newcommand{\m}{\mu}              

\newcommand{\n}{\nu}              \newcommand{\cP}{\mathcal{P}}

\renewcommand{\r}{\rho}             

\newcommand{\s}{\sigma}           \newcommand{\cR}{\mathcal{R}}

           \newcommand{\cS}{\mathcal{S}}

\renewcommand{\t}{\tau}             

\newcommand{\f}{\phi}             

\newcommand{\F}{\Phi}             

\newcommand{\vp}{\varphi}

\newcommand{\p}{\psi}             
 
             \newcommand{\cZ}{\mathcal{Z}}
      
\renewcommand{\o}{\omega}

\newcommand{\vk}{\varkappa}

  \def\mA{{\mathscr A}}
 \def\mB{{\mathscr B}}

  \def\mH{{\mathscr H}}

  \def\mP{{\mathscr P}}

  \def\mU{{\mathscr U}}

\newcommand{\gD}{\mathfrak{D}}

\newcommand{\gR}{\mathfrak{R}}
\newcommand{\gS}{\mathfrak{S}}

\def\Z{\mathbb{Z}}
\def\R{\mathbb{R}}
\def\C{\mathbb{C}}

\def\N{\mathbb{N}}

\def\qqq{\qquad}
\def\qq{\quad}

\let\ge\geqslant
\let\le\leqslant

\newcommand{\ca}{\begin{cases}}
\newcommand{\ac}{\end{cases}}
\newcommand{\ma}{\begin{pmatrix}}
\newcommand{\am}{\end{pmatrix}}

\def\lt{\biggl}
\def\rt{\biggr}

\renewcommand{\[}{\begin{equation}}
\renewcommand{\]}{\end{equation}}

\def\wt{\widetilde}
\def\pa{\partial}
\def\sm{\setminus}
\def\es{\emptyset}
\def\no{\noindent}
\def\ol{\overline}
\def\iy{\infty}
\def\ev{\equiv}
\def\/{\over}
\def\ts{\times}

\def\os{\oplus}
\def\ss{\subset}

\def\Re{\mathop{\rm Re}\nolimits}

\def\Im{\mathop{\rm Im}\nolimits}
\def\supp{\mathop{\rm supp}\nolimits}
\def\sign{\mathop{\rm sign}\nolimits}

\def\Tr{\mathop{\rm Tr}\nolimits}
\def\const{\mathop{\rm const}\nolimits}

\def\BBox{\hspace{1mm}\vrule height6pt width5.5pt depth0pt \hspace{6pt}}
\def\wh{\widehat}

\begin{document}
\title{Magnetic Schr\"odinger operators on
armchair nanotubes}
\author{
Andrey Badanin
\begin{footnote} {
Arkhangelsk State Technical University,
e-mail: a.badanin@agtu.ru }
\end{footnote}
\and
Evgeny Korotyaev
\begin{footnote} {corresponding author,
Institut f\"ur  Mathematik,  Humboldt Universit\"at zu Berlin,
e-mail: evgeny@math.hu-berlin.de\ \ }
\end{footnote}
}

\maketitle

\begin{abstract}
\no We consider the Schr\"odinger operator with a periodic potential  on
a quasi 1D continuous periodic  model of armchair nanotubes in $\R^3$ in
a uniform magnetic field (with amplitude $B\in \R$), which  is parallel
to the axis of the  nanotube.
The spectrum of this operator consists of an absolutely continuous part
(spectral bands separated by gaps) plus an infinite number of
eigenvalues  with infinite multiplicity.
We describe all eigenfunctions with the same eigenvalue including
compactly supported.
We describe the spectrum as a function of $B$.
For some specific potentials
we prove an existence of gaps independent on the magnetic field.
If $B\ne 0$, then there exists
an infinite number of gaps $G_n$ with the length $|G_n|\to\iy$
as $n\to\iy$, and we determine the asymptotics of the gaps
at high energy for fixed $B$. Moreover, we determine
the asymptotics of the gaps $G_n$ as $B\to 0$ for fixed $n$.

\end{abstract}

\section{Introduction  and main results}
\setcounter{equation}{0}

\begin{figure}
\centering
\noindent
\tiny
\hfill
\includegraphics[height=.4\textheight]{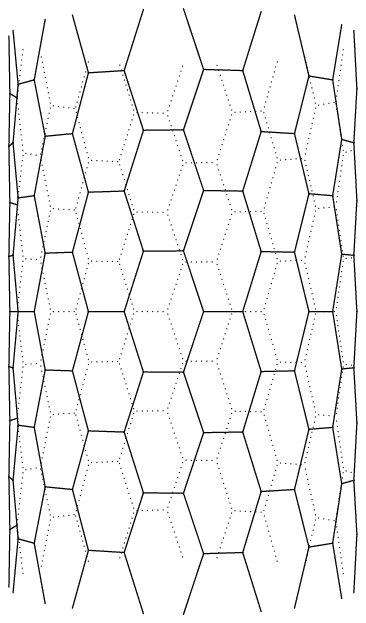}
\hfill
{
\psfrag{g010}{$\Gamma_{0,1}$}
\psfrag{g020}{$\Gamma_{0,2}$}
\psfrag{g030}{$\Gamma_{0,3}$}
\psfrag{g040}{$\Gamma_{0,4}$}
\psfrag{g050}{$\Gamma_{0,5}$}
\psfrag{g060}{$\Gamma_{0,6}$}
\psfrag{g110}{$\Gamma_{1,1}$}
\psfrag{g120}{$\Gamma_{1,2}$}
\psfrag{g130}{$\Gamma_{1,3}$}
\psfrag{g140}{$\Gamma_{1,4}$}
\psfrag{g150}{$\Gamma_{1,5}$}
\psfrag{g160}{$\Gamma_{1,6}$}
\psfrag{g-110}{$\Gamma_{-1,1}$}
\psfrag{g-120}{$\Gamma_{-1,2}$}
\psfrag{g-130}{$\Gamma_{-1,3}$}
\psfrag{g-140}{$\Gamma_{-1,4}$}
\psfrag{g-150}{$\Gamma_{-1,5}$}
\psfrag{g-160}{$\Gamma_{-1,6}$}
\includegraphics[height=.4\textheight]{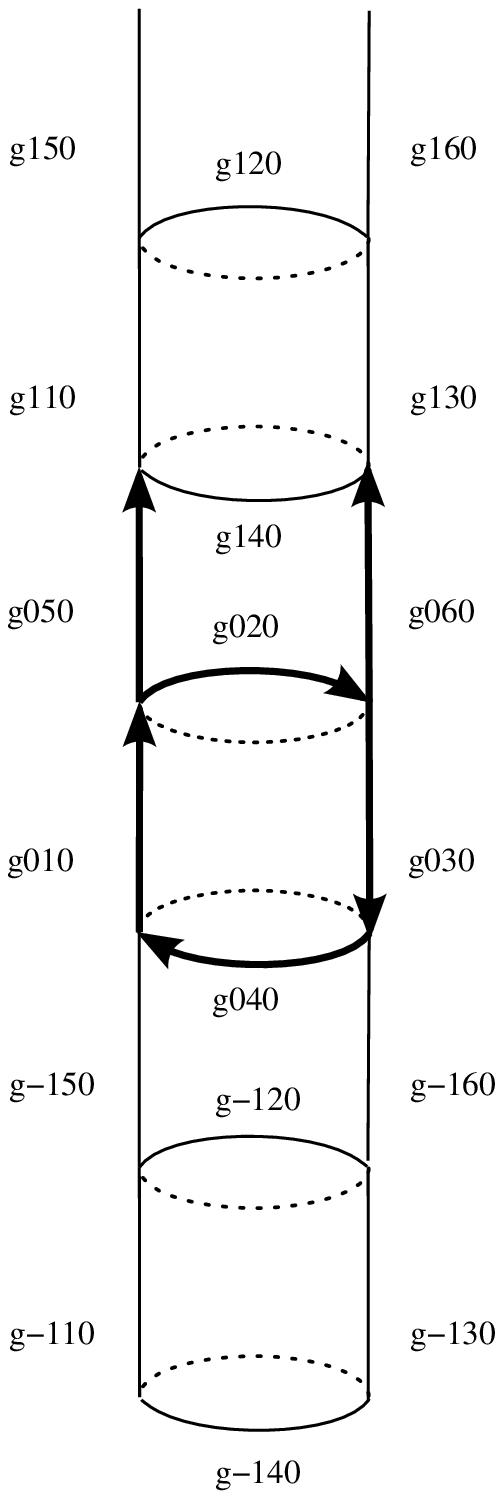}
}
\hfill~
\caption{Armchair graph for $N=10$ and for $N=1$.}
\lb{fig1}
\end{figure}

We consider the Schr\"odinger operator
$\mH_B =(-i\nabla-\mA)^2+V_q$  with a periodic potential $V_q$
on the armchair nanotube $\G^N\ss\R^3,N\ge 1$
in a uniform magnetic field
$\mB=B(0,0,1)\in \R^3$, $B\in\R$.
The corresponding vector potential is given by
$\mA(x)={1\/2}[\mB,x]={b\/2}(-x_2,x_1,0),\  x=(x_1,x_2,x_3)\in\R^3$.
Our model nanotube $\G^N$ is a union of edges $\G_\o$, i.e.,
$$
\G^N=\cup_{\o\in \cZ} \G_\o,\qq \o=(n,j,k)\in \cZ=\Z\ts \N_6\ts \Z_N,
\qq \N_m=\{1,2,...,m\},\qq \Z_N=\Z/(N\Z),
$$
  see Fig. \ref{fig1}, \ref{fig2}. Each edge $\G_\o=\{x=
{\bf r}_\o+t{\bf e}_\o,\  t\in [0,1]\}$  is oriented
by the vector ${\bf  e}_\o\in \R^3$ and has starting point
$\bf r_\o\in \R^3$.
We have the coordinate $x={\bf r}_\o+t{\bf e}_\o$ and the local
coordinate $t\in [0,1]$ (length preserving).
We define ${\bf r}_{\o}, {\bf e}_{\o},\o=(n,j,k)\in \cZ$ by
$$
{\bf r}_{n,1,k}=2nh{\bf e}_3+R(c_{2k}, s_{2k},0),\qq
{\bf r}_{n,4,k}=2nh{\bf e}_3+R(\cos \b_k, \sin \b_k,0),\
 \b_k=2\b+\f_{2k},
$$
$$
{\bf r}_{n,2,k}={\bf r}_{n,5,k}=(2n+1)h{\bf e}_3+R(\cos\z_k,\sin\z_k,0),\qq
\z_k=\b-\a+\f_{2k},\qq \f_k={\pi k\/N},
$$$$
{\bf r}_{n,3,k}={\bf r}_{n,6,k}=(2n+1)h{\bf e}_3+R(c_{2k+1}, s_{2k+1},0),\
$$
$$
{\bf e}_{n,j,k}={\bf r}_{n,j+1,k}-{\bf r}_{n,j,k},\qq j=1,2,3,\qq
{\bf e}_{n,4,k}={\bf r}_{n,1,k+1}-{\bf r}_{n,4,k},\
{\bf e}_{n,5,k}={\bf r}_{n+1,1,k}-{\bf r}_{n,5,k},\qq
$$
\[
\lb{ab}
{\bf e}_{n,6,k}={\bf r}_{n+1,4,k}-{\bf r}_{n,6,k},\ \
 \sin \b={1\/R},\ \ \sin \a={1\/2R},\ \
R={\sqrt{\cos\f_1+{5\/4}}\/\sin\f_1},
\]
$0<\a<\b, \a+\b={\pi\/N}$ and $h=\sqrt{2+R_1R_2-2R^2}$,
$R_p=\sqrt{(pR)^2-1}, p=1,2$.

Each edge $\G_\o$ is the segment  and connect the
vertices, which lie on the cylinder
$\cC\ev\{(x_1,x_2,x_3)\in \R^3:x_1^2+x_2^2=R^2\}$.
Note that each edge $\G_\o$ (without the endpoints) lies inside the
cylinder $\cC$.
The points ${\bf r}_{0,j,k}$ are vertices of the regular N-gon
$\mP_j, 1=1,2,3,4$.
$\mP_3$ ($\mP_2$ ) arises from $\mP_1$  ($\mP_4$ )
by the following motion:
rotate around the axis of the cylinder $\cC$ by the angle ${\pi\/N}$
and translate by $h{\bf e}_3$.
Repeating this procedure we obtain $\G^N$.

For each function $y$ on $\G^N$
we define a function $y_\o=y|_{\G_\o}, \o\in \cZ$.
We identify each function $y_\o$
on $\G_\o$ with a function on $[0,1]$ by using the local coordinate
$t\in [0,1]$.
Define the Hilbert space $L^2(\G^N)=\os_{\o\in \cZ} L^2(\G_\o)$.
Our operator $\mH_B$ on $\G^N$ acts in the Hilbert space
$L^2(\G^N)=\os_\o  L^2(\G_\o)$ and is given by
\[
(\mH_B f)_\o=-\pa_{\o}^2f_\o(t)+q(t)f_\o(t),\qq
\pa_{\o}={d\/dt}-ia_\o,
\qq a_\o(t)=(\mA({\bf r}_\o+t{\bf e}_\o),{\bf e}_\o),
\]
see \cite{SDD},
where  $(V_q f)_\o=qf_\o, q\in L^2(0,1)$  and
$\os_\o f_\o, \os_\o f_\o'' \in L^2(\G^N)$ satisfies

\no {\bf Magnetic Kirchhoff Boundary Conditions:} {\it $y$ is
continuous on $\G^N$ and
\[
\lb{KirC}
\pa_{\o_2}y_{\o_2}(0)-\pa_{\o_1}y_{\o_1}(1)+\pa_{\o_5}y_{\o_5}(0)=0,\qqq
\pa_{\o_3}y_{\o_3}(0)-\pa_{\o_2}y_{\o_2}(1)+\pa_{\o_6}y_{\o_6}(0)=0,
$$
$$
-\pa_{\o_3}y_{\o_3}(1)+\pa_{\o_4}y_{\o_4}(0)-\pa y_{n-1,6,k}(1)=0,\qqq
-\pa_{\o_4}y_{\o_4}(1)+\pa y_{n,1,k+1}(0)-\pa  y_{n-1,5,k+1}(1)=0,
\]
for all $ \o_j=(n,j,k), (n,j,k)\in \Z\ts\N_6\ts\Z_N$}.

Condition \er{KirC} means  that  the sum of derivatives of $y$
at each vertex of $\G^N$ equals 0 and the orientation of edges gives
the sign $\pm$.
Such models  were introduced by Pauling \cite{Pa}
in 1936 to simulate aromatic molecules. They were described
in more detail by Ruedenberg and Scherr \cite{RS} in 1953.
Further progress had been made toward periodic systems
by Coulson in \cite{Co} where a network model of graphite
layer was worked out.
A network model of a crystal with non-trivial
potentials along bonds was studied by Montroll \cite{Mo}.
Carbon and boron nitride nanostructures,
in particular nanotubes, graphene layers have been modelled
by quantum networks with a honeycomb lattice structure
in \cite{ALM}, \cite{ALM1}.
As it is shown in \cite{Ha}, Ch.7.3.3, the structure of
carbon-boron-nitride single wall nanotubes may be complicated.
The edge potential $q$ is not even for such tubes.
We consider one of nanotube models, where the potential is
not even. The other models can be considered by the similar way.

\begin{figure}
\centering
\noindent
{
\tiny
\psfrag{g011}{$\Gamma_{0,1,1}$}
\psfrag{g021}{$\Gamma_{0,2,1}$}
\psfrag{g031}{$\Gamma_{0,3,1}$}
\psfrag{g041}{$\Gamma_{0,4,1}$}
\psfrag{g051}{$\Gamma_{0,5,1}$}
\psfrag{g061}{$\Gamma_{0,6,1}$}
\psfrag{g012}{$\Gamma_{0,1,2}$}
\psfrag{g022}{$\Gamma_{0,2,2}$}
\psfrag{g032}{$\Gamma_{0,3,2}$}
\psfrag{g042}{$\Gamma_{0,4,2}$}
\psfrag{g052}{$\Gamma_{0,5,2}$}
\psfrag{g062}{$\Gamma_{0,6,2}$}
\psfrag{g01N}{$\Gamma_{0,1,N}$}
\psfrag{g02N}{$\Gamma_{0,2,N}$}
\psfrag{g03N}{$\Gamma_{0,3,N}$}
\psfrag{g04N}{$\Gamma_{0,4,N}$}
\psfrag{g05N}{$\Gamma_{0,5,N}$}
\psfrag{g06N}{$\Gamma_{0,6,N}$}

\psfrag{g111}{$\Gamma_{1,1,1}$}
\psfrag{g121}{$\Gamma_{1,2,1}$}
\psfrag{g131}{$\Gamma_{1,3,1}$}
\psfrag{g141}{$\Gamma_{1,4,1}$}
\psfrag{g151}{$\Gamma_{1,5,1}$}
\psfrag{g161}{$\Gamma_{1,6,1}$}
\psfrag{g112}{$\Gamma_{1,1,2}$}
\psfrag{g122}{$\Gamma_{1,2,2}$}
\psfrag{g132}{$\Gamma_{1,3,2}$}
\psfrag{g142}{$\Gamma_{1,4,2}$}
\psfrag{g152}{$\Gamma_{1,5,2}$}
\psfrag{g162}{$\Gamma_{1,6,2}$}
\psfrag{g11N}{$\Gamma_{1,1,N}$}
\psfrag{g2N}{$\Gamma_{1,2,N}$}
\psfrag{g13N}{$\Gamma_{1,3,N}$}
\psfrag{g14N}{$\Gamma_{1,4,N}$}
\psfrag{g15N}{$\Gamma_{1,5,N}$}
\psfrag{g16N}{$\Gamma_{1,6,N}$}
\psfrag{g-111}{$\Gamma_{-1,1,1}$}
\psfrag{g-121}{$\Gamma_{-1,2,1}$}
\psfrag{g-131}{$\Gamma_{-1,3,1}$}
\psfrag{g-141}{$\Gamma_{-1,4,1}$}
\psfrag{g-151}{$\Gamma_{-1,5,1}$}
\psfrag{g-161}{$\Gamma_{-1,6,1}$}
\psfrag{g-112}{$\Gamma_{-1,1,2}$}
\psfrag{g-122}{$\Gamma_{-1,2,2}$}
\psfrag{g-132}{$\Gamma_{-1,3,2}$}
\psfrag{g-142}{$\Gamma_{-1,4,2}$}
\psfrag{g-152}{$\Gamma_{-1,5,2}$}
\psfrag{g-162}{$\Gamma_{-1,6,2}$}
\psfrag{g-11N}{$\Gamma_{-1,1,N}$}
\psfrag{g-12N}{$\Gamma_{-1,2,N}$}
\psfrag{g-13N}{$\Gamma_{-1,3,N}$}
\psfrag{g-14N}{$\Gamma_{-1,4,N}$}
\psfrag{g-15N}{$\Gamma_{-1,5,N}$}
\psfrag{g-16N}{$\Gamma_{-1,6,N}$}
\psfrag{O1}{$\Omega_{1}$}
\psfrag{O2}{$\Omega_{2}$}
\includegraphics[width=.8\textwidth,height=.5\textwidth]{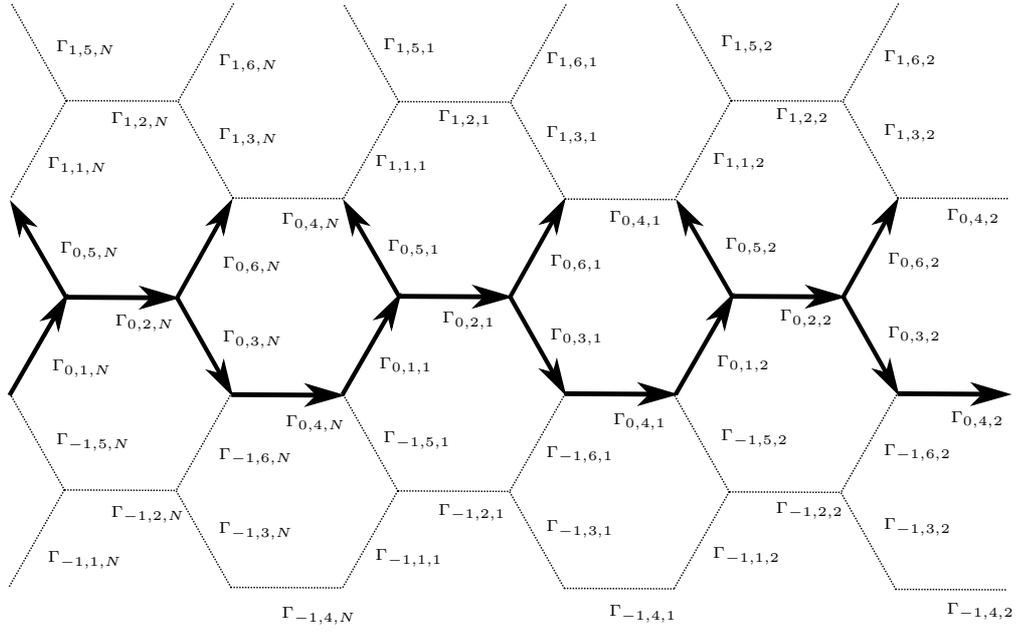}
}
\caption{ A piece of a nanotube $\G^N$.
The fundamental domain is marked by a bold line.
}
\lb{fig2}
\end{figure}

The standard arguments (see \cite{KL}) yield  that $\mH_B$ is
self-adjoint.
For simplicity we shall denote $\G_{\a,1}\ss \G^1$ by  $\G_{\a}$, for
$\a=(n,j)\in \cZ_1=\Z\ts \N_6$. Thus $\G^1=\cup_{\a\in \cZ_1} \G_\a$,
see Fig \ref{fig1}. In Theorem \ref{T1} we will show that
$\mH_B$ is unitarily equivalent to $H^a=\os_1^N H_k^a$,
where the operator $H_k^a$ acts in the Hilbert space $L^2(\G^1)$
and is given by
$(H_k^a f)_\a=-f_\a''+q f_\a,\ (f_\a)_{\a\in \cZ_1},
(f_\a'')_{\a\in \cZ_1}\in L^2(\G^1)$,
and the components $f_\a$  satisfy:
\begin{multline}
\label{1K0}
e^{ia_1}f_{n,1}(1)=f_{n,2}(0)=f_{n,5}(0),\qq
e^{ia_2}f_{n,2}(1)=f_{n,3}(0)=f_{n,6}(0),\\
e^{ia_1}f_{n,3}(1)=f_{n,4}(0)=e^{ia_1}f_{n-1,6}(1),\qq
e^{ia_2}s^k f_{n,4}(1)=f_{n,1}(0)=e^{-ia_1}f_{n-1,5}(1),
\end{multline}
\begin{multline}
\label{1K1}
e^{ia_1}f_{n,1}'(1)-f_{n,2}'(0)-f_{n,5}'(0)=0,\qq
e^{ia_2}f_{n,2}'(1)-f_{n,3}'(0)-f_{n,6}'(0)=0,\\
e^{ia_1}f_{n,3}'(1)-f_{n,4}'(0)+e^{ia_1}f_{n-1,6}'(1)=0,\qq
e^{ia_2}s^k f_{n,4}'(1)-f_{n,1}'(0)+e^{-ia_1}f_{n-1,5}'(1)=0,
\end{multline}
where
\[
\lb{a12}
s=e^{2\pi i\/N},\qq
a_1={B(R_2-R_1)\/4},
\qq
a_2={BR_2\/4},\qq
a=a_1+a_2={B(2R_2-R_1)\/4},
\]
$R$ is given by \er{ab}, and $R_p=\sqrt{(pR)^2-1},p=1,2$.
For the operator $H_k^a$ we define fundamental solutions
$\vt_k^{(\n)}=(\vt_{k,\a}^{(\n)})_{\a\in \cZ_1}$,
$\vp_k^{(\n)}=(\vp_{k,\a}^{(\n)})_{\a\in \cZ_1}$,
$\n=1,2$  which satisfy
\begin{multline}
\lb{eqf0}
-\D y+V_qy=\l y, \qq  on \ \ \G^1,\qqq
\text{b. c. \er{1K0},\er{1K1}}\\
\vT_{k,0}(1,\l)=\F_{k,0}'(1,\l)=I_2,
\qq \vT_{k,0}'(1,\l)=\F_{k,0}(1,\l)=0,
\end{multline}
where $I_n,n\ge 2$ is the $n\ts n$ identity matrix and
\[
\lb{dmm1}
\vT_{k,n}=\ma\vt_{k,n,5}^{(1)} & \vt_{k,n,5}^{(2)} \\
\vt_{k,n,6}^{(1)} & \vt_{k,n,6}^{(2)}\am,\qq
\F_{k,n}=\ma\vp_{k,n,5}^{(1)} & \vp_{k,n,5}^{(2)} \\
\vp_{k,n,6}^{(1)} & \vp_{k,n,6}^{(2)}\am.
\]
We define the monodromy matrix by
\[
\lb{dmm}
\cM_k(\l)=\ma \vT_{k,1} & \F_{k,1}\\
\vT_{k,1}' & \F_{k,1}'\am(1,\l).
\]
This matrix is determined by the fundamental
solutions on a fundamental cell $\G_0^1=\cup_{j=1}^6 \G_{0,j}$.
In order to analyse the operator $H^a$ we use the methods based on
a research of the monodromy matrix $\cM_k$ (see \cite{KL},
\cite{KL1},\cite{BBKL},\cite{BBK1},
and also \cite{BBK}, \cite{CK}, \cite{K}).

Let
$\s_D=\{\m_n, n\ge 1\}$ be the spectrum of the problem $-y''+qy=\l y,
y(0)=y(1)=0$ (the Dirichlet spectrum).
We formulate our first result.

\begin{theorem}
\label{T1}
(i) The operator $\mH_B$ is unitarily equivalent to
$H^a=\os_1^N H_k^a$.

\no (ii) Let $(a,k)\in\R\ts\Z_N$.
Then for any $\l\in\C\sm\s_D$
there exist unique fundamental
solutions $\vt_k^{(\n)}=(\vt_{k,\a}^{(\n)})_{\a\in \cZ_1},
\vp_k^{(\n)}=(\vp_{k,\a}^{(\n)})_{\a\in \cZ_1},\n=1,2$.
 Moreover, each of the functions $\vt_{k,\a}^{(\n)}(x,\l)$,
$\vp_{k,\a}^{(\n)}(x,\l),x\in \G^1$ is
analytic in $\l\in\C\sm\s_D$ and the monodromy matrix
$\cM_k(\l),\l\not\in\s_D$
has eigenvalues of the form $\t_{k,1}^{\pm 1 }, \t_{k,2}^{\pm 1 }$
and satisfies
\[
\lb{T1-3}
\det \cM_k=1,\qqq
\cM_k(\l)^\top  J\cM_k(\l)=J,\qq where \qq J=\ma 0& {\bf j}\\
-{\bf j}& 0\am,\ \ {\bf j}=\ma 0& 1\\ 1& 0\am.
\]
Furthermore, the matrix-valued function
$\cR\cM_k\cR^{-1}$ is entire, where $\cR=I_2\os \vp_1I_2$.

\end{theorem}

\no For all $B\in\R$ the monodromy matrix $\cM_k$ has a simple pole
at each point $\l\in \s_D$, which is an eigenvalue of $H_k^a$
(see Theorem \ref{T2}) and has no any other singularities.
Moreover, these poles are independent on the magnetic field.
This result is in contrast with the corresponding result
for the zigzag tube (see \cite{KL}, \cite{KL1}).
For the zigzag tube there exist some singular values of the
magnetic field, where the monodromy matrix is not well defined.
Moreover, the spectrum is pure point for
such singular values and the spectral bands shrink to the points
as the magnetic field approach to the singular value.
%

%
%

For the equation $-y''+q(t)y=\l y$ on the real line we introduce the
fundamental solutions $\vt(t,\l)$ and $\vp(t,\l),(t,\l)\in\R\ts\C$
satisfying $\vt(0,\l)=\vp'(0,\l)=1, \vt'(0,\l)=\vp(0,\l)=0$. We
define the functions $F$, $F_-$ by
\[
\lb{Ho}
F={\vp_1'+\vt_1\/2},\qqq  F_-={\vp_1'-\vt_1\/2},
\]
where
$\vp_1=\vp(1,\cdot),\vt_1=\vt(1,\cdot),\vp_1'=\vp'(1,\cdot),
\vt_1'=\vt'(1,\cdot)$.
We formulate the results about the Lyapunov function and the spectrum
of $H_k^a$.

\begin{theorem}
\lb{T3}  Let $(a,k)\in\R\ts\Z_N$
and let $\t_{k,1}^{\pm 1 }, \t_{k,2}^{\pm 1 }$
be eigenvalues of $\cM_k$. Then

\no (i) The Lyapunov functions
$F_{k,\nu}={1\/2}(\t_{k,\nu}+\t_{k,\nu}^{-1}), \n=1,2$
are branches of $F_k=\xi_k+\sqrt{\r_k}$
on the two sheeted Riemann surface $\gR_k$ defined by $\sqrt {\r_k}$
and satisfy:
\[
\lb{DeLk}
F_{k,\nu}=\xi_k-(-1)^\n\sqrt{\r_k},\qq
\xi_k={9F^2-F_-^2-1\/2}-s_k^2,\qq
\r_k=(9F^2-s_k^2)c_k^2+s_k^2F_-^2,
\]
\[
\lb{S3b}
\det(\cM_k-\t I_4)
=(\t^2-2F_{k,1}\t+1)(\t^2-2F_{k,2}\t+1),
\]
where
\[
\lb{cs}
c_k=\cos({\pi k\/N}+a),\qqq s_k=\sin({\pi k\/N}+a),
\]
$a$ is given by \er{a12}.

\no (ii) If $F_k(\l)\in (-1,1)$ for some $\l\in \R$ and
$\l$ is not a branch point of $F_k(\l)$, then $F_k'(\l)\neq 0$.

\no (iii)  The  following identities hold:
\begin{multline}
\lb{T3-1}
\s(H_k^a)=\s_{pp}(H_k^a)\cup\s_{ac}(H_k^a),\qq
\s_{pp}(H_k^a)=\s_D,\\
\s_{ac}(H_k^a)=\{\l\in\R: F_{k,\nu}(\l)\in [-1,1]\ \text{for some}\
\nu\in\N_2\},
\end{multline}
\[
\lb{sh1}
\s(H_k^{-a})=\s(H_{N-k}^a),\ \
\s(H_k^{a+{\pi\/N}})=\s(H_{k+1}^a),\ \
\s(H_k^{a+\pi})=\s(H_{k}^a)
\]
(each point of the spectrum is counted with multiplicities).
Each point $\l_0\in \s_{pp}(H_k^a)$ has
an infinite multiplicity.

\end{theorem}

\no {\bf Remark.} 1) We take the branch of $\sqrt{\r_k}$ such that
if $\r_k(\l)>0,\l\in\R$, then $\sqrt{\r_k(\l)}>0$ and
$F_{k,1}=\xi_k+\sqrt{\r_k}>
F_{k,2}=\xi_k-\sqrt{\r_k}$.
Then $\sqrt{\r_k}$ is not an entire function even for $c_k=0$ or $s_k=0$.
Using other branches
we could define a new function $\sqrt{\r_k}$ which is entire
for $c_k=0$ or $s_k=0$, but this choice is not convenient
for labeling of endpoints of gaps (see below Theorems \ref{Tk}-\ref{Te}).

\no 2) The statement (ii) is similar to the standard result
of Lyapunov for the Hamilton systems.
The  similar results hold for the zigzag tube in
magnetic field (see \cite{KL},\cite{KL1}).

\no 3) Periodicity of the spectrum with respect
to $B$ holds for the zigzag nanotube
(see \cite{KL1}).

\no 4) Below we will sometimes write
$c_k(a),F_{k,\n}(\l,a),...$
instead of $c_k,F_{k,\n}(\l),...$,
when several magnetic fields are being dealth with.
In order to prove \er{sh1} we use the identities
\[
\lb{scs}
s_k(-a)=s_{N-k}(a),\ \ c_k(-a)=c_{N-k}(a),\ \
s_k(a+{\pi\/N})=s_{k+1}(a),\ \ c_k(a+{\pi\/N})=c_{k+1}(a).
\]

\no 5) We say that $\{\l_0\}\ss\s_{pp}(H_k^a)$ is a flat band.
The corresponding eigenfunctions of $H_k^a$ have compact supports
(shown by Fig.\ref{sef}) and are described
in Theorem \ref{T2}.

We define the entire functions
\[
\lb{S3n}
D_k^\pm=\det(\cM_k\mp I_4)=4(F_{k,1}\mp 1)(F_{k,2}\mp 1).
\]
\no The zeros $\l_{\n,2n}^{k,\pm}, n\ge 0,\n=1,2$, of the function
$D_k^+$ are the periodic eigenvalues.
The zeros $\l_{\n,2n-1}^{k,\pm},(\n,n)\in\N_2\ts\N$, of
$D_k^-$ are the antiperiodic eigenvalues.
Let
$$
\l_{2,0}^{k,+}\le\l_{1,0}^{k,+}
\le\l_{1,2}^{k,-}\le\l_{2,2}^{k,-}\le
\l_{2,2}^{k,+}\le\l_{1,2}^{k,+}\le
\l_{1,4}^{k,-}\le\l_{2,4}^{k,-}\le
\l_{2,4}^{k,+}\le\l_{1,4}^{k,+}\le...
$$
and
$$
\l_{2,1}^{k,-}\le\l_{1,1}^{k,-}\le\l_{1,1}^{k,+}\le\l_{2,1}^{k,+}
\le\l_{2,3}^{k,-}\le\l_{1,3}^{k,-}\le\l_{1,3}^{k,+}\le\l_{2,3}^{k,+}\le...
$$
counted with multiplicities.
Such labeling was introduced in \cite{BBK1} for the case $a=0$,
and is associated with
the Lyapunov functions $F_{k,1},F_{k,2}$ (see Fig.\ref{lfd}).

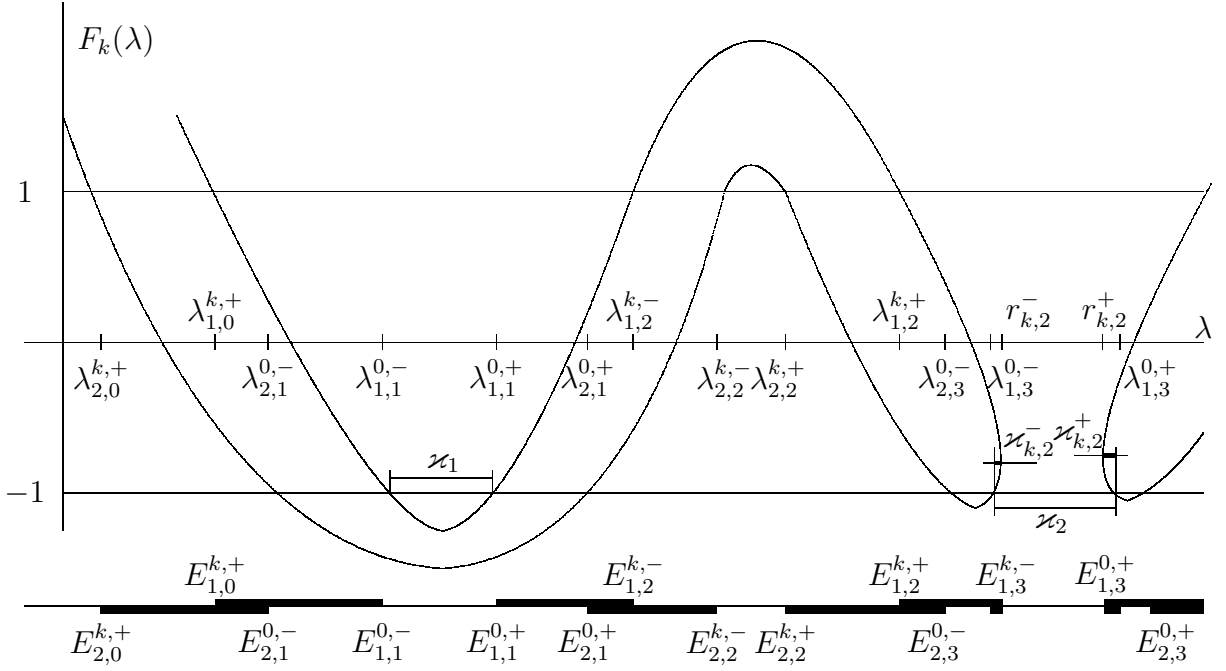
\begin{figure}
\unitlength 1.00mm \linethickness{0.2pt}
\begin{picture}(162.00,80.00)(00.00,-20.00)
\put(10.00,-5.00){\line(0,1){70.00}}
\put(10.00,40.00){\line(1,0){150.00}}
\put(05.00,20.00){\line(1,0){155.00}}
\put(10.00,00.00){\line(1,0){150.00}}
\put(160.00,22.00){\makebox(0,0)[cc]{$\l$}}
\put(05.00,00.00){\makebox(0,0)[cc]{$-1$}}
\put(05.00,40.00){\makebox(0,0)[cc]{$1$}}
\put(17.00,60.00){\makebox(0,0)[cc]{$F_k(\l)$}}
\bezier{600}(10.00,50.00)(30.00,-8.00)(60.00,-10.00)
\bezier{600}(60.00,-5.00)(50.00,-4.00)(25.00,50.00)
\bezier{400}(60.00,-5.00)(70.00,-3.00)(85.00,40.00)
\bezier{400}(85.00,40.00)(100.00,80.00)(120.00,40.00)
\bezier{400}(120.00,40.00)(140.00,02.00)(130.00,-2.00)
\bezier{400}(60.00,-10.00)(85.00,-8.00)(97.00,40.00)
\bezier{150}(97.00,40.00)(100.00,47.00)(105.00,40.00)
\bezier{400}(130.00,-2.00)(120.00,02.00)(105.00,40.00)
\bezier{400}(150.00,-1.00)(140.00,02.00)(161.00,41.00)
\bezier{250}(150.00,-1.00)(155.00,01.00)(160.00,08.00)
\put(133.50,19.00){\line(0,1){2.00}}
\put(137.00,24.00){\makebox(0,0)[cc]{$r_{k,2}^-$}}
\put(146.70,19.00){\line(0,1){2.00}}
\put(147.00,24.00){\makebox(0,0)[cc]{$r_{k,2}^+$}}
\put(15.00,19.00){\line(0,1){2.00}}
\put(15.00,15.00){\makebox(0,0)[cc]{$\l_{2,0}^{k,+}$}}
\put(30.00,19.00){\line(0,1){2.00}}
\put(30.00,24.00){\makebox(0,0)[cc]{$\l_{1,0}^{k,+}$}}
\put(85.00,19.00){\line(0,1){2.00}}
\put(85.00,24.00){\makebox(0,0)[cc]{$\l_{1,2}^{k,-}$}}
\put(120.00,19.00){\line(0,1){2.00}}
\put(120.00,24.00){\makebox(0,0)[cc]{$\l_{1,2}^{k,+}$}}
\put(96.00,19.00){\line(0,1){2.00}}
\put(97.00,15.00){\makebox(0,0)[cc]{$\l_{2,2}^{k,-}$}}
\put(105.00,19.00){\line(0,1){2.00}}
\put(104.00,15.00){\makebox(0,0)[cc]{$\l_{2,2}^{k,+}$}}
\put(37.00,19.00){\line(0,1){2.00}}
\put(37.00,15.00){\makebox(0,0)[cc]{$\l_{2,1}^{0,-}$}}
\put(52.00,19.00){\line(0,1){2.00}}
\put(52.00,15.00){\makebox(0,0)[cc]{$\l_{1,1}^{0,-}$}}
\put(67.00,19.00){\line(0,1){2.00}}
\put(67.00,15.00){\makebox(0,0)[cc]{$\l_{1,1}^{0,+}$}}
\put(79.00,19.00){\line(0,1){2.00}}
\put(79.00,15.00){\makebox(0,0)[cc]{$\l_{2,1}^{0,+}$}}
\put(126.00,19.00){\line(0,1){2.00}}
\put(126.00,15.00){\makebox(0,0)[cc]{$\l_{2,3}^{0,-}$}}
\put(132.00,19.00){\line(0,1){2.00}}
\put(135.00,15.00){\makebox(0,0)[cc]{$\l_{1,3}^{0,-}$}}
\put(149.00,19.00){\line(0,1){2.00}}
\put(153.00,15.00){\makebox(0,0)[cc]{$\l_{1,3}^{0,+}$}}
\put(53.00,00.00){\line(0,1){3.00}}
\put(66.50,00.00){\line(0,1){3.00}}
\put(60.00,04.00){\makebox(0,0)[cc]{$\vk_1$}}
\put(53.00,02.00){\line(1,0){13.50}}
\put(132.50,-03.00){\line(0,1){9.00}}
\put(148.50,-03.00){\line(0,1){9.00}}
\put(132.50,-02.00){\line(1,0){16.00}}
\put(140.00,-04.00){\makebox(0,0)[cc]{$\vk_2$}}
\put(131.00,04.00){\line(1,0){7.00}}
\put(132.50,04.20){\line(1,0){0.70}}
\put(132.50,03.80){\line(1,0){0.70}}
\put(137.00,07.00){\makebox(0,0)[cc]{$\vk_{k,2}^-$}}
\put(143.00,05.00){\line(1,0){7.00}}
\put(147.00,05.20){\line(1,0){1.50}}
\put(147.00,04.80){\line(1,0){1.50}}
\put(143.50,08.00){\makebox(0,0)[cc]{$\vk_{k,2}^+$}}
\put(05.00,-15.00){\line(1,0){155.00}}
\put(15.00,-15.90){\line(1,0){22.00}}
\put(15.00,-15.80){\line(1,0){22.00}}
\put(15.00,-15.70){\line(1,0){22.00}}
\put(15.00,-15.60){\line(1,0){22.00}}
\put(15.00,-15.50){\line(1,0){22.00}}
\put(15.00,-15.40){\line(1,0){22.00}}
\put(15.00,-15.30){\line(1,0){22.00}}
\put(15.00,-15.20){\line(1,0){22.00}}
\put(30.00,-14.80){\line(1,0){22.00}}
\put(30.00,-14.70){\line(1,0){22.00}}
\put(30.00,-14.60){\line(1,0){22.00}}
\put(30.00,-14.50){\line(1,0){22.00}}
\put(30.00,-14.40){\line(1,0){22.00}}
\put(30.00,-14.30){\line(1,0){22.00}}
\put(30.00,-14.20){\line(1,0){22.00}}
\put(30.00,-14.10){\line(1,0){22.00}}
\put(67.00,-14.10){\line(1,0){18.00}}
\put(67.00,-14.20){\line(1,0){18.00}}
\put(67.00,-14.30){\line(1,0){18.00}}
\put(67.00,-14.40){\line(1,0){18.00}}
\put(67.00,-14.50){\line(1,0){18.00}}
\put(67.00,-14.60){\line(1,0){18.00}}
\put(67.00,-14.70){\line(1,0){18.00}}
\put(67.00,-14.80){\line(1,0){18.00}}
\put(79.00,-15.90){\line(1,0){17.00}}
\put(79.00,-15.80){\line(1,0){17.00}}
\put(79.00,-15.70){\line(1,0){17.00}}
\put(79.00,-15.60){\line(1,0){17.00}}
\put(79.00,-15.50){\line(1,0){17.00}}
\put(79.00,-15.40){\line(1,0){17.00}}
\put(79.00,-15.30){\line(1,0){17.00}}
\put(79.00,-15.20){\line(1,0){17.00}}
\put(105.00,-15.90){\line(1,0){21.00}}
\put(105.00,-15.80){\line(1,0){21.00}}
\put(105.00,-15.70){\line(1,0){21.00}}
\put(105.00,-15.60){\line(1,0){21.00}}
\put(105.00,-15.50){\line(1,0){21.00}}
\put(105.00,-15.40){\line(1,0){21.00}}
\put(105.00,-15.30){\line(1,0){21.00}}
\put(105.00,-15.20){\line(1,0){21.00}}
\put(120.00,-14.10){\line(1,0){13.50}}
\put(120.00,-14.20){\line(1,0){13.50}}
\put(120.00,-14.30){\line(1,0){13.50}}
\put(120.00,-14.40){\line(1,0){13.50}}
\put(120.00,-14.50){\line(1,0){13.50}}
\put(120.00,-14.60){\line(1,0){13.50}}
\put(120.00,-14.70){\line(1,0){13.50}}
\put(120.00,-14.80){\line(1,0){13.50}}
\put(132.00,-15.90){\line(1,0){1.50}}
\put(132.00,-15.80){\line(1,0){1.50}}
\put(132.00,-15.70){\line(1,0){1.50}}
\put(132.00,-15.60){\line(1,0){1.50}}
\put(132.00,-15.50){\line(1,0){1.50}}
\put(132.00,-15.40){\line(1,0){1.50}}
\put(132.00,-15.30){\line(1,0){1.50}}
\put(132.00,-15.20){\line(1,0){1.50}}
\put(147.00,-14.10){\line(1,0){13.00}}
\put(147.00,-14.20){\line(1,0){13.00}}
\put(147.00,-14.30){\line(1,0){13.00}}
\put(147.00,-14.40){\line(1,0){13.00}}
\put(147.00,-14.50){\line(1,0){13.00}}
\put(147.00,-14.60){\line(1,0){13.00}}
\put(147.00,-14.70){\line(1,0){13.00}}
\put(147.00,-14.80){\line(1,0){13.00}}
\put(147.00,-15.90){\line(1,0){2.00}}
\put(147.00,-15.80){\line(1,0){2.00}}
\put(147.00,-15.70){\line(1,0){2.00}}
\put(147.00,-15.60){\line(1,0){2.00}}
\put(147.00,-15.50){\line(1,0){2.00}}
\put(147.00,-15.40){\line(1,0){2.00}}
\put(147.00,-15.30){\line(1,0){2.00}}
\put(147.00,-15.20){\line(1,0){2.00}}
\put(153.00,-15.90){\line(1,0){7.00}}
\put(153.00,-15.80){\line(1,0){7.00}}
\put(153.00,-15.60){\line(1,0){7.00}}
\put(153.00,-15.50){\line(1,0){7.00}}
\put(153.00,-15.40){\line(1,0){7.00}}
\put(153.00,-15.30){\line(1,0){7.00}}
\put(153.00,-15.20){\line(1,0){7.00}}
\put(153.00,-15.70){\line(1,0){7.00}}
\put(134.00,-11.00){\makebox(0,0)[cc]{$E_{1,3}^{k,-}$}}
\put(147.00,-11.00){\makebox(0,0)[cc]{$E_{1,3}^{0,+}$}}
\put(15.00,-20.00){\makebox(0,0)[cc]{$E_{2,0}^{k,+}$}}
\put(30.00,-11.00){\makebox(0,0)[cc]{$E_{1,0}^{k,+}$}}
\put(85.00,-11.00){\makebox(0,0)[cc]{$E_{1,2}^{k,-}$}}
\put(120.00,-11.00){\makebox(0,0)[cc]{$E_{1,2}^{k,+}$}}
\put(96.00,-20.00){\makebox(0,0)[cc]{$E_{2,2}^{k,-}$}}
\put(105.00,-20.00){\makebox(0,0)[cc]{$E_{2,2}^{k,+}$}}
\put(37.00,-20.00){\makebox(0,0)[cc]{$E_{2,1}^{0,-}$}}
\put(52.00,-20.00){\makebox(0,0)[cc]{$E_{1,1}^{0,-}$}}
\put(67.00,-20.00){\makebox(0,0)[cc]{$E_{1,1}^{0,+}$}}
\put(79.00,-20.00){\makebox(0,0)[cc]{$E_{2,1}^{0,+}$}}
\put(125.00,-20.00){\makebox(0,0)[cc]{$E_{2,3}^{0,-}$}}
\put(155.00,-20.00){\makebox(0,0)[cc]{$E_{2,3}^{0,+}$}}
\end{picture}
\caption{Graph of the function $F_k(\l)$ and the spectrum of $H_k$}
\lb{lfd}
\end{figure}

A zero of $\r_k, k\in \Z_N$ is called a {\bf resonance} of $H_k^a$.
There exist real and non-real resonances
(see example from \cite{BBKL}). Roughly speaking the simple real
resonances create gaps. Note that in the case of
zigzag nanotubes all resonances are real
(see \cite{KL}, \cite{KL1}).

We define the functions
\[
\lb{uv}
u_k(\l)=|F_-(\l)|-s_k^2,\qq v_k(\l)=|F_-(\l)|-c_k^2,\qq
\l\in\R,\qq k\in\Z_N.
\]

\begin{theorem}
\lb{Tk}  Let $(a,k)\in\R\ts\Z_N$. Then
the identity
\[
\lb{sHk}
\s_{ac}(H_k^a)=\cup_{(\n,n)\in\N_2\ts\N}S_{\n,n}^{k,a},\qq
S_{\n,n}^{k,a}=[E_{\n,n-1}^{k,+},E_{\n,n}^{k,-}]
\]
holds, where each spectral band $S_{\n,n}^{k,a}$
satisfies:
\[
\lb{s}
E_{\n,p-1}^{k,\pm}=\l_{\n,p-1}^{k,\pm},\qq
E_{2,p}^{k,\pm}=\l_{2,p}^{0,\pm},\qq
E_{1,p}^{k,\pm}=\ca
\l_{1,p}^{0,\pm}\qq \text{if}\qq
v_k(\l_{1,p}^{0,\pm})\ge 0\\
r_{k,n}^\pm\qq\ \text{if}\qq
v_k(\l_{1,p}^{0,\pm})<0
\ac\!\!\!,\qq p=2n-1,
\]
\[
\lb{gc}
\l_{\n,p}^{k,\pm}(a)=\l_{\n,p}^{0,\pm}(0),\qq
E_{2,p}^{k,\pm}(a)=E_{2,p}^{0,\pm}(0),\qq
\vk_n(a)=\vk_n(0),
\]
where
\[
\lb{kan}
r_{k,n}^-=\min\{\l\in\ol\vk_n:\r_k(\l)=0\},\ \
r_{k,n}^+=\max\{\l\in\ol\vk_n:\r_k(\l)=0\},\ \
\vk_n=(\l_{1,p}^{0,-},\l_{1,p}^{0,+}).
\]

\end{theorem}

\no {\bf Remark.}
1) For some specific potentials and some $n\ge 1$
the intervals $\vk_n$ are antiperiodic gaps
(see below Proposition \ref{pro}
and Theorem \ref{Te} (i)). These gaps are independent on
the magnetic field.

\no 2) We describe $\s_{ac}(H_k^a)$ in details
in Theorems \ref{4s}, \ref{C}.

Recall the needed properties of the Hill operator $\wt H
y=-y''+q(t)y$ on the real line with a periodic potential
$q(t+1)=q(t),t\in \R$. The spectrum of $\wt H$ is purely absolutely
continuous and consists of intervals
$[\wt\l_{n-1}^+,\wt\l_n^-], n\ge 1$. These intervals are
separated by the gaps $\wt\g_n=(\wt\l_n^-,\wt\l_n^+)$ of length
$|\wt\g_n|\ge 0$. If a gap $\wt\g_n$ is degenerate, i.e.
$|\wt\g_n|=0$, then the corresponding segments merge.
The sequence $\wt\l_0^+<\wt\l_1^-\le \wt\l_1^+\ <...$ is the
spectrum of the equation $-y''+qy=\l y$ with 2-periodic boundary
conditions,  that is  $y(t+2)=y(t), t\in \R$.
Here equality $\wt\l_n^-= \wt\l_n^+$ means that $\wt\l_n^\pm$ is an
eigenvalue of
multiplicity 2.  It is well-known that
$\m_n\in [\wt\l^-_n,\wt\l^+_n ],n\ge 1$.
The function $F$ has only simple zeros $\e_n,n\ge 1$, which satisfy
$\e_1<\e_2<...$
We describe the spectral gaps of $H_k^a$.

\begin{theorem}\lb{g} Let $(a,k)\in\R\ts\Z_N$.
Then

\no (i) $\s_{ac}(H_k^a)=\R\sm\cup_{n\ge 0} G_{k,n}^a$, where the gaps
$G_{k,n}^a$ satisfy:
\begin{multline}
\lb{Hk1}
\wt\g_0\ss G_{k,0}^a=(-\iy,E_{2,0}^{k,+}),\ \
\wt\g_n\ss G_{k,4n}^a=(E_{2,2n}^{k,-},E_{2,2n}^{k,+}),\ \
G_{k,4n-2}^a=(E_{1,2n-1}^{k,-},E_{1,2n-1}^{k,+})\ss\vk_n,\\
G_{k,4n-3}^a=(E_{2,2n-1}^{k,-},E_{1,2n-2}^{k,+}),\qq
G_{k,4n-1}^a=(E_{1,2n}^{k,-},E_{2,2n-1}^{k,+}),\qq
\e_n\in[E_{1,2n-1}^{k,-},E_{1,2n-1}^{k,+}].
\end{multline}
If $s_k^2(a)<s_\ell^2(a_1)$
for some $(a_1,\ell)\in\R\ts\Z_N$, then
\[
\lb{gk2}
G_{k,4n}^a\ss G_{\ell,4n}^{a_1},\qqq
G_{\ell,2n-1}^{a_1}\ss G_{k,2n-1}^a\qqq\text{all}\qqq n\ge 1.
\]
Moreover, for all $n\ge n_0$ and some $n_0\ge 1$

\[
\lb{gk4}
G_{k,4n-2}^a=\ca \es,\qq\ \text{if}\ (a,k)=(0,0)\\
\ne\es,\ \text{if}\ (a,k)\ne(0,0)
\ac\!\!\!,\qq
G_{k,2n-1}^a=\es,\ \text{if}\ (a,k)\ne (0,0).
\]

\no (ii) If $v_k=|F_-(\l)|-c_k^2\ge 0$ on $\vk_n$ for some $n\ge 1$, then
$G_{k,4n-2}^a=\vk_n$.

\end{theorem}

We formulate our main result about the spectrum of $H^a$.

\begin{theorem} \lb{TM}
Let $a\in[0,{\pi\/{2N}}]$. Then
the following identities hold true:
\[
\lb{TM-1}
\s(H^a)=\s_{pp}(H^a)\cup\s_{ac}(H^a),\qq
\s_{pp}(H^a)=\s_D,\qq \s(H^a)=\s(H^{-a})=\s(H^{a+{\pi\/N}})
\]
(each point of the spectrum is counted with multiplicities).

\no (i) $\s_{ac}(H^a)=\R\sm\cup_{n\ge 0} G_{n}^a$,
where the gaps
$G_n^a=\cap_{k\in\Z_N}G_{k,n}^a=(E_{n}^-,E_{n}^+)$ satisfy:
\[
\lb{TM-3}
\wt\g_{n}\ss G_{4n}^a,
\qq G_{4n-2}^a\ss\vk_n,\qq
\e_n\in [E_{4n-2}^-,E_{4n-2}^+],
\]
\[
\lb{g3}
G_{4n}^a\ss G_{4n}^{a_1},\qq G_{2n-1}^{a_1}\ss G_{2n-1}^a,\qqq
0\le a<a_1\le {\pi\/2N}.
\]
Moreover, for all $n\ge n_0$ and some $n_0\ge 1$
\[
\lb{pgH}
G_{4n-2}^a=\ca \es ,\qq\  \text{if}\ a=0\\
\ne\es , \  \text{if}\  a\ne 0
\ac\!\!\!,\qqq
G_{2n-1}^a=\es,\ \text{if}\ (a,N)\ne (0,1).
\]

\no (ii) The spectrum $\s(H^a),a\in(0,{\pi\/2N})$, has
an infinite number of gaps $G_{2n}^a\ne\es$
and $|G_{2n}^a|\to\iy$ as $n\to\iy$.
Moreover, the following asymptotics hold true:
\[
\lb{aE}
E_{4n-2m}^\pm(a)=\lt(\pi (n-{m\/2})\pm \wt\theta_m\rt)^2+q_0+O(n^{-1}),
\qq n\to\iy,\qq m=0,1,
\]
where $q_0=\int_0^1q(t)dt,\wt\theta_0=\arccos({1\/3}\sqrt{5+4|\cos a|}),
\wt\theta_1=\arcsin({1\/3}\sin a)$.
Furthemore, if $N\ge 2$, then all other gaps
with sufficiently large energy are empty.

\no (iii) If  $|F_-|\ge\cos^2a$
on $\vk_n$ for some $n\ge 1$, then
$G_{4n-2}^a=\vk_n$.

\end{theorem}

\no {\bf Remark.} 1) Identities \er{TM-1} show that for analysis
of operator $H^a,a\in\R$ it is sufficiently to study $a\in[0,{\pi\/2N}]$.

\no 2) Using the statement (iii) we prove in Proposition \ref{pro}
that each $G_{4n-2}^a=G_{4n-2}^0\ne\es$, $(a,n)\in\R\ts\N_{n_0}$
for any $n_0\ge 1$ and some non-even potential (which depends on $n_0$).
Note that if the potential is even, then all the gaps $G_{4n-2}^a$
depend on the magnetic field (see below \er{rge}).

Recall the definition of the effective masses for the Hill operator.
The identity $F(\l)=\cos k(\l), \l\in \C_+$ defines the quasimomentum
$k(\l)$, which is an analytic function on $\l\in \C_+$.
With  each endpoint  of  the gap $\wt\g_n\ne \es$, we associate
the effective mass $M_0^+, M_n^{\pm}$ by
(we take some branches $k$ such that $k(\l)\to \pi n$
as $\l\to \wt\l_n^{\pm}$)
\[
\l=\wt\l_n^{\pm}+{(k(\l)-\pi n)^2\/2M_n^{\pm}}(1 +o(1))
\ \ \ \ \  {\rm as} \ \ \ \             \l\to \wt\l_n^{\pm}.
\]
(see \cite{KK}).
We formulate results for even potentials.
Let $L_{even}^2(0,1)=\{q\in L^2(0,1):q(1-t)=q(t),t\in(0,1)\}$.

\begin{theorem} \lb{Te} Let $q\in L_{even}^2(0,1)$,
$(a,n)\in[0,{\pi\/2N}]\ts\N$. Then

\no (i) If $c_k=0$, then $\r_k=0$ and $G_{k,4n-2}^a=\vk_n$.
If $c_k\ne 0$, then all resonances are real and
$G_{k,4n-2}^a=(r_{k,n}^-,r_{k,n}^+)$.
 For each $k\in\Z_N$ the gaps $G_{k,4n-2}^a$ satisfy:
\[
\lb{gke}
G_{k,4n-2}^a=\ca \es,\qq\ \text{if}\qq (a,k)=(0,0)\\
\ne\es,\qq \text{if}\qq (a,k)\ne(0,0)\!\!\!
\ac,\qqq G_{k,2n-1}^a=\es.
\]
If $s_k^2(a)<s_\ell^2(a_1)$
for some $(a_1,\ell)\in\R\ts\Z_N$, then
$G_{k,4n-2}^a\ss G_{\ell,4n-2}^{a_1}.$

\no (ii) The gaps of $H^a$ satisfy
\[
\lb{mfg}
G_{4n-2}^a\ss G_{4n-2}^{a_1},\qq
0\le a<a_1\le{\pi\/2N},
\]
\[
\lb{rge}
G_{4n}^0=\wt\g_n,\qqq G_{4n-2}^a=\ca \es ,\qq  \text{if}\ a=0\\
\ne\es , \  \text{if}\  a\ne 0\ac\!\!\!,\qqq
G_{2n-1}^a=\es.
\]
Moreover, if $a\in(0,{\pi\/2N}]$, then $G_{2n}^a\ne\es$, and
\[
\lb{as1}
E_{4n}^\pm(a)=\wt\l_n^\pm+{a^2\/9M_n^\pm}+O(a^4),
\qq \text{if}\qq \wt\g_n\ne\es,
\]
\[
\lb{as2}
E_{4n}^\pm(a)=\wt\l_n^+
\pm {\sqrt 2 a\/3\sqrt{|F''(\wt\l_n^+)|}}+O(a^2),
\qq \text{if}\qq \wt\g_n=\es,
\]
\[
\lb{as3}
E_{4n-2}^\pm(a)=\e_{n}\pm {a\/3|F'(\e_{n})|}+O(a^2)
\]
as $a\to 0$.

\end{theorem}

Schr\"odinger operator
on graphs is a subject of intensive studies in the last years
(see \cite{BGP}, \cite{KL}, \cite{KS}, \cite{KuP},...).
In \cite{BGP} the magnetic Schr\"odinger
operator on the planar square lattice was considered.
The spectrum of operator is expressed in terms of the Lyapunov function
of the corresponding Hill operator. In \cite{P} these results
are extended to the case of an arbitrary connected $3D$
graph with identical edges and the same even potential.

Schr\"odinger operators on nanotubes have attracted
a lot of attention recently
(see \cite{KL}, \cite{KL1}, \cite{KuP},...).
Authors of \cite{KuP} consider the case
of the zigzag, armchair and chiral nanotubes with even
potentials $q$. They show that the spectrum of the operator, as a set,
coincides with the spectrum of the Hill operator. Note that these results
can be obtained from the results of \cite{P}
(K.Pankrashkin, non-published).
In some sence the case of even potentials $q$ is closed to the case
$q=0$. Indeed, if $q$ is even, then $F_-=0$ \cite{MW} and
by Theorem \ref{T3}, the corresponding Lyapunov functions are
expressed in terms of the Lyapunov function $F$ for the Hill operator.
Using the identity $F(\l)=\cos\k(\l)$, where $\k(\l)$ is
a quasimomentum for the Hill operator, the Lyapunov function
is expressed in terms of $\cos\k(\l)$,
which is similar to the case $q=0$, where $F(\l)=\cos\sqrt\l$.

The detailed descripiton of spectrum of Schr\"odinger operator and
magnetic Schr\"odinger operator on zigzag nanotube with arbitrary
potentials was given in the articles \cite{KL}, \cite{KL1}, \cite{K1}.
The operator was represented as the direct sum of quasi
one-dimensional operators. Research of the spectrum
of these quasi one-dimensional operators is based on
the analysis of the monodromy matrix (see \cite{BBK}, \cite{CK}).
Authors of \cite{KL} describe all eigenfunctions with the same
eigenvalue. They define a Lyapunov function,
which is analytic on some Riemann surface. On each
sheet, the Lyapunov function has the same properties
as in the scalar case, but it has
branch points (resonances).
They prove that all  resonances are real and they determine the
asymptotics of the periodic and anti-periodic spectrum and of the
resonances at high energy. They show that there exist two types of gaps:
i) stable gaps, where the endpoints are periodic and
anti-periodic eigenvalues, ii) unstable (resonance) gaps, where the
endpoints are resonances (i.e., real branch points of the Lyapunov
function).
They describe all finite gap potentials. They show that the mapping:
 potential $\to$ all eigenvalues is a real analytic isomorphism
for some class of potentials.

In \cite{KL1} authors describe how the spectrum depends
on the magnetic field. They observe so-called localization
of the spectrum for some values of the magnetic field, when
the spectrum of the operator become pure point. This phenominon
is typical for quantum graphs and was discussed in many papers
(see ref. in \cite{P1}).

Korotyaev \cite{K1} considers the effective masses
for zigzag nanotubes in magnetic fields.
Following \cite{KL}, \cite{KL1},
he consider so called modified Lyapunov
function and the corresponding quasimomentum.
The modified Lyapunov function is entire,
then the quasimomentum is an analytic function in $\C_+$.
As for the case of Hill operator,
the quasimomentum is a conformal mapping of $\C_+$
onto the domain $\Re k>0,\Im k>0$ with vertical slits.
Identities and a priori estimates
of gap lengths, heights of slits  in terms of effective masses
are obtained. These values as functions of magnetic field are
described.

Schr\"odinger operator on armchair nanotube with arbitrary potentials
was considered in \cite{BBKL}, \cite{BBK1}. In \cite{BBKL}
all eigenfunctions of the operator with the same
eigenvalue are described. The Lyapunov function is analytic on some
Riemann surface and has branch points (resonances).
Some example is considered in this paper,
which shows that there are real and non-real resonances.
The detailed analysis of the spectrum is a subject of the article
\cite{BBK1}. There exist different types of gaps:
i) periodic (antiperiodic) gaps, where the endpoints are periodic
(antiperiodic) eigenvalues,
ii) resonance gaps, where the endpoints are resonances,
iii) r-mix gaps, where one endpoint is an antiperiodic eigenvalue
and other endpoint is a resonance,
iv) p-mix gaps, where one endpoint is an antiperiodic eigenvalue
and other endpoint is a periodic eigenvalue.
Asymptotics of the gaps at high energy are obtained.

We present the plan of the paper. In Section 2 we prove
Theorems \ref{T1}, \ref{T3} and Theorem \ref{T2} about the eigenfunctions.
In Section 3 we prove the main Theorems \ref{Tk}-\ref{Te}.
Furthermore, in this Section we prove Theorems \ref{4s}, \ref{C},
where we give the detail analysis of the spectrum including
the multiplicity, description of endpoints of gaps etc.
Moreover, in this Section we prove an existence of gaps independent on the magnetic field for some specific potentials
(see Proposition \ref{pro}).

\section{Spectrum of the operators $H_k^a$}
\setcounter{equation}{0}

We need the following result.

\begin{lemma}
\lb{mafi}
Each  function
$a_\o(t)=(\mA({\bf r}_\o+t{\bf e}_\o),{\bf e}_\o)=a_j,\
\o=(n,j,k)\in \cZ,\  t\in [0,1]$, where
\[
\lb{aoI}
a_1=a_{3}=-a_{5}=a_{6}={bR^2\/2}\sin (\b-\a)={b(R_2-R_1)\/4},
\qq a_2=a_{4}={bR^2\/2}\sin 2\a={bR_2\/4}.
\]

\end{lemma}

\no{\bf Proof.}
Identity $\mA({\bf r})={b\/2}[{\bf e}_3,{\bf r}], {\bf e}_3=(0,0,1),
{\bf r}\in \R^3$ yields
\[
\label{ao2}
a_\o(t)={b\/2}([{\bf e}_3,{\bf r}_\o+t{\bf e}_\o],{\bf e}_\o)=
{b\/2}([{\bf e}_3,{\bf r}_\o],{\bf e}_\o)=a_\o(0)=a_\o,\qqq any
\qq t\in[0,1].
\]
Moreover, the vector field $\mA({\bf r})$ is parallel to the plane $x_3=0$
(this plane is orthogonal to the axis of the cylinder $\cC$) and $\mA({\bf r})$ is tangential
to the cylinder surface $x_1^2+x_2^2=R^2$. We deduce that each projection
$a_\o,\o\in \cZ$ of $\mA({\bf r})$ to the edge $\G_\o$ satisfies $a_\o=a_j$
and $a_1=a_{3}=-a_{5}=a_{6},a_2=a_{4}$.

We have ${\bf e}_{0,j,k}={\bf r}_{0,j+1,k}-{\bf r}_{0,j,k}, j=1,2$, which yields
$$
a_{0,j,k}={b\/2}([{\bf e}_3,{\bf r}_{0,j,k}],{\bf r}_{0,j+1,k}).
$$
If $j=1$, then ${\bf r}_{0,1,k}=R(c_{2k},s_{2k},0), {\bf r}_{0,2,k}=h{\bf e}_3
+R(\cos\z_k,\sin\z_k,0)$ and thus
$$
a_1=a_{0,1,k}={bR^2\/2}\det
\ma \cos\z_k & \sin\z_k & 0\\
0 & 0 & 1\\
c_{2k} & s_{2k} & 0\am
={bR^2\/2}\sin (\z_k-\f_{2k})={bR^2\/2}\sin (\b-\a).
$$
If $j=2$, then ${\bf r}_{0,3,k}=h{\bf e}_3+R(c_{2k+1},s_{2k+1},0)$ and thus
$$
a_2=a_{0,2,k}={bR^2\/2}\det
\ma c_{2k+1}& s_{2k+1} & 0\\
0 & 0 & 1\\
\cos\z_k & \sin\z_k & 0\am={bR^2\/2}\sin (\f_{2k+1}-\z_k)
={bR^2\/2}\sin 2\a,
$$
since $\f_{2k+1}-\z_k={\pi\/N}-\b+\a=2\a$, which yields \er{aoI}.
$\BBox$

\no Now we will prove Theorem \ref{T1} and the identities
\[
\lb{T1-2}
\Tr\cM_0
=2(9F^2-F_-^2-1),\quad
\Tr\cM_k=\Tr\cM_0-4s_k^2,
\]
\[
\lb{T1-4}
\Tr\cM_0^2
=72F^2+{1\/2}(\Tr\cM_0)^2-4,\qq
\Tr\cM_k^2=\Tr\cM_0^2-8s_k^2\Tr\cM_0-16s_k^2c_k^2.
\]

\no {\bf Proof of Theorem \ref{T1} and identities \er{T1-2}, \er{T1-4}.}
(i) Introduce the unitary operator $\mU$ in $L^2(\G^N)$
by $(\mU f)_\o=e^{ita_\o}f_\o, f=(f_\o)_{\o\in \cZ}$, where
$a_\o$ is given by \er{aoI}.
We define a modified operator $H^a=\mU^* \mH_B \mU$ in $L^2(\G^N)$.
Then   $H^a$ is  given by
$(H^af)_\o=-f_\o''+qf_\o, \o\in \cZ$, where
$f\in \gD(H^a)$. The domain $\gD(H^a)$ consists of the functions $f =
(f_\o)_{\o\in\cZ}, (f_\o'')_{\o\in\cZ}\in L^2(\G^N)$ and  $f$ satisfies
{\bf the Modified Kirchhoff  Boundary Conditions:}
\begin{multline}
\label{KirC0}
e^{ia_1}f_{\o_1}(1)=f_{\o_2}(0)=f_{\o_5}(0),\quad
e^{ia_2}f_{\o_2}(1)=f_{\o_3}(0)=f_{\o_6}(0),\\
e^{ia_1}f_{\o_3}(1)=f_{\o_4}(0)=e^{ia_1}f_{n-1,6,k}(1),\quad
e^{ia_2}f_{\o_4}(1)=f_{n,1,k+1}(0)=e^{-ia_1}f_{n-1,5,k+1}(1),
\end{multline}
\begin{multline}
\lb{KirC1}
f'_{\o_2}(0)-e^{ia_1}f'_{\o_1}(1)+f'_{\o_5}(0)=0,\qq
f'_{\o_3}(0)-e^{ia_2}f'_{\o_2}(1)+f'_{\o_6}(0)=0,\\
-e^{ia_1}f'_{\o_3}(1)+f'_{\o_4}(0)-e^{ia_1}f'_{n-1,6,k}(1)=0,\qq
-e^{ia_2}f'_{\o_4}(1)+f'_{n,1,k+1}(0)-e^{-ia_1}f'_{n-1,5,k+1}(1)=0,
\end{multline}
for all $ \o_j=(n,j,k), j=\N_6, (n,k)\in \Z\ts\Z_N$.

Define the operator $\cS$ in $\C^N$ by
$\cS u=(u_N,u_1,\dots,u_{N-1})^\top$, $u=(u_n)_1^N\in \C^N$.
The unitary operator $\cS$ has the form $\cS =\sum_1^Ns^k\cP_k$,
where $\cS e_k=s^ke_k$,
$e_k={1\/N^{1\/2}}(1,s^{-k},s^{-2k},...,s^{-kN+k})$
 is an eigenvector (recall $s=e^{i{2\pi \/N}}$);
 $\cP_ku=e_k(u,e_k)$ is a projector. The function $f$
 in the Kirchhoff boundary conditions \er{KirC0}
is a vector function  $f=(f_{\o}), \o=(n,j,k)\in\cZ$.  We define a new
vector-valued function $f_{n,j}=(f_{n,j,k})_{k=1}^{N}\in \C^N, (n,j)\in
\cZ_1=\Z\ts \N_6$, which  satisfies the equation
$-f_{n,j}''+qf_{n,j}=\l f_{n,j}$,
and the conditions
\begin{multline}
\label{C0}
e^{ia_1}f_{n,1}(1)=f_{n,2}(0)=f_{n,5}(0),\quad
e^{ia_2}f_{n,2}(1)=f_{n,3}(0)=f_{n,6}(0),\\
e^{ia_1}f_{n,3}(1)=f_{n,4}(0)=e^{ia_1}f_{n-1,6}(1),\quad
e^{ia_2}\cS f_{n,4}(1)=f_{n,1}(0)=e^{-ia_1}f_{n-1,5}(1),
\end{multline}
\begin{multline}
\label{C1}
e^{ia_1}f_{n,1}'(1)-f_{n,2}'(0)-f_{n,5}'(0)=0,\qq
e^{ia_2}f_{n,2}'(1)-f_{n,3}'(0)-f_{n,6}'(0)=0,\\
e^{ia_1}f_{n,3}'(1)-f_{n,4}'(0)+e^{ia_1}f_{n-1,6}'(1)=0,\qq
e^{ia_2}\cS f_{n,4}'(1)-f_{n,1}'(0)+e^{-ia_1}f_{n-1,5}'(1)=0,
\end{multline}
for all $n\in Z$, which follow from the Kirchhoff conditions \er{KirC}.
The operators $\cS$ and $\mH_B$ commute, then  we deduce that
$\mH_B\cP_k$ is unitarily equivalent to the operator $H_k^a$ acting
 in $L^2(\G^1)$ and $H_k^a$ is given by
$ (H_k^af)_\a=-f_\a''+q(t)f_\a$, where $(f_\a)_{\a\in \cZ_1},
(f_\a'')_{\a\in \cZ_1}\in L^2(\G^1)$ and components $f_\a$
satisfy the boundary conditions \er{1K0}, \er{1K1}.
Thus $\mH_B$ is unitarily equivalent to the operator $H^a=\os_1^N H_k^a$.

\no Proof of (ii) and \er{T1-2}, \er{T1-4}
repeats the arguments from \cite{BBKL}.
$\BBox$

Define the subspace $\cH_k^a(\l)=\{\p\in \gD(H_k^a):
H_k^a\p=\l \p\}$
for $(\l,k)\in \s_{pp}(H_k^a)\ts\Z_N$.

\begin{theorem} \lb{T2}
Let $(a,\l,k)\in\R\ts\s_D\ts \Z_N$. Then

\no i) Every eigenfunction of $\cH_k^a(\l)$ vanishes at all vertices
of $\G^1$.

\no ii)
There exist the functions
$\p^{(0,\n)}=(\p^{(0,\n)}_\a)_{\a\in\cZ_1}\in\cH_k^a(\l),\n\in\N_2$ on $\G^1$ such that
$\supp \p^{(0,\n)}$ $\ss\cup_{\ell=0,1}(\cup_{j\in\N_6}\G_{\ell,j})$,
each function
$\p^{(n,\n)}=(\p^{(n,\n)}_{m,j})_{(m,j)\in \cZ_1}
=(\p^{(0,\n)}_{m-n,j})_{(m,j)\in \cZ_1}\in \cH_k^a(\l),
n\in \Z$. Moreover, each $f\in \cH_k^a(\l)$ satisfies:
\[
\lb{T2-3}
f=\sum_{(n,\n)\in\Z\ts\N_2} \wh  f_{n,\n}\p^{(n,\n)},\qqq
(\wh f_{n,1},\wh f_{n,2})_{n\in\Z}\in \ell^2\os\ell^2,
\]
where if $(s_k,\f^2)=(0,1),$ $\f=\vp_1'(\l)$, then
\[
\lb{fna}
\wh f_{n,1}=f_{n,1}'(0),\qqq
\wh  f_{n,2}=f_{n,2}'(0),
\]
if $(s_k,\f^2)\ne (0,1)$, then
\[
\lb{T2-4}
\wh f_{n,1}={f_{n,5}'(0)+s^k\f e^{i(a+a_1)}f_{n,6}'(0)\/\wt\vk_1\wt\vk_2},\qq
\wh f_{n,2}=-{e^{ia_1}f_{n,6}'(0)+\f e^{ia}f_{n,5}'(0)\/\wt\vk_1\wt\vk_2},
\]
$\wt\vk_1=1-s^ke^{2ia}\f^2\ne 0$, $\wt\vk_2=1-s^ke^{2ia}\f^4\ne 0$.
Moreover, the mapping $f\to (\wh f_{n,1},\wh f_{n,2})_{n\in\Z}$
is a linear isomorphism between $\cH_k^a(\l)$ and $\ell^2\os\ell^2$.


\end{theorem}

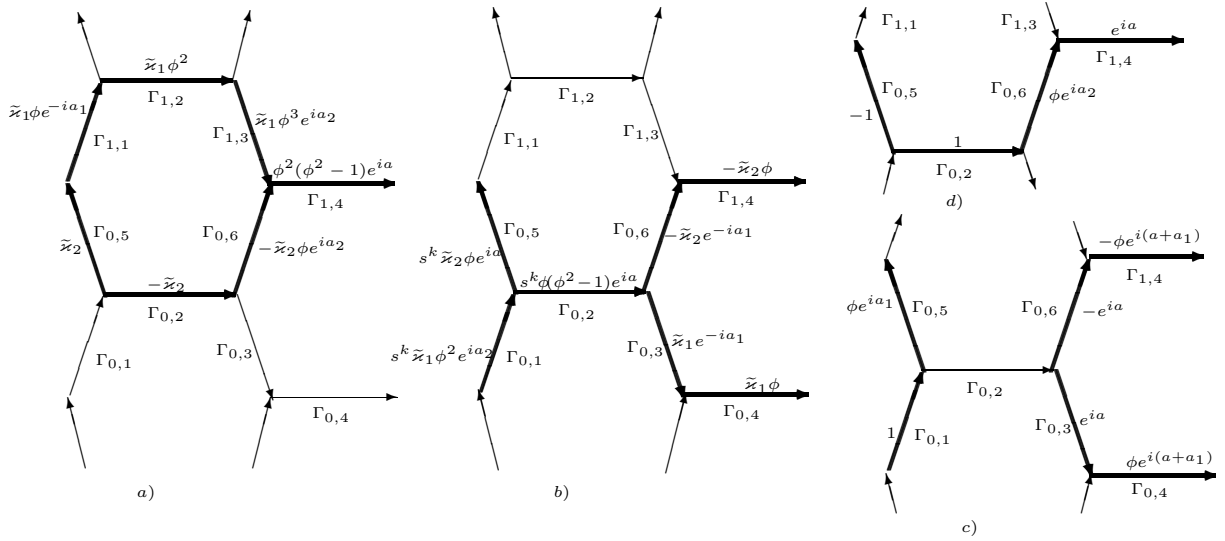
\begin{figure}
\unitlength=1.00mm
\special{em:linewidth 0.4pt}
\linethickness{0.2pt}
\tiny
\begin{picture}(159.01,68.67)
\put(8.32,16.67){\vector(1,3){4.33}}
\put(12.66,30.00){\vector(1,0){17.33}}
\put(12.66,30.10){\vector(1,0){17.33}}
\put(12.66,29.90){\vector(1,0){17.33}}
\put(12.66,30.20){\vector(1,0){17.33}}
\put(12.66,29.80){\vector(1,0){17.33}}
\put(20.99,30.33){\makebox(0,0)[cb]{$-\wt\vk_2$}}
\put(34.66,16.34){\vector(1,0){16.67}}
\put(29.82,29.67){\vector(1,3){5.00}}
\put(29.72,29.67){\vector(1,3){5.00}}
\put(29.92,29.67){\vector(1,3){5.00}}
\put(29.62,29.67){\vector(1,3){5.00}}
\put(30.02,29.67){\vector(1,3){5.00}}
\put(32.99,36.33){\makebox(0,0)[lc]{$-\wt\vk_2\f e^{ia_2}$}}
\put(12.99,29.67){\vector(-1,3){5.00}}
\put(13.09,29.67){\vector(-1,3){5.00}}
\put(12.89,29.67){\vector(-1,3){5.00}}
\put(13.19,29.67){\vector(-1,3){5.00}}
\put(12.79,29.67){\vector(-1,3){5.00}}
\put(9.99,36.33){\makebox(0,0)[rc]{$\wt\vk_2$}}
\put(7.99,45.00){\vector(1,3){4.33}}
\put(8.09,45.00){\vector(1,3){4.33}}
\put(7.89,45.00){\vector(1,3){4.33}}
\put(8.19,45.00){\vector(1,3){4.33}}
\put(7.79,45.00){\vector(1,3){4.33}}
\put(10.99,54.33){\makebox(0,0)[rc]{$\wt\vk_1\!\f e^{\!-ia_{1}}$}}
\put(12.32,58.34){\vector(1,0){17.33}}
\put(12.32,58.24){\vector(1,0){17.33}}
\put(12.32,58.44){\vector(1,0){17.33}}
\put(12.32,58.14){\vector(1,0){17.33}}
\put(12.32,58.54){\vector(1,0){17.33}}
\put(20.99,59.33){\makebox(0,0)[cb]{$\wt\vk_1\f^2$}}
\put(34.32,44.67){\vector(1,0){16.67}}
\put(34.32,44.57){\vector(1,0){16.67}}
\put(34.32,44.77){\vector(1,0){16.67}}
\put(34.32,44.47){\vector(1,0){16.67}}
\put(34.32,44.87){\vector(1,0){16.67}}
\put(42.99,45.33){\makebox(0,0)[cb]{$\f^2(\f^2-1)e^{ia}$}}
\put(29.99,58.33){\vector(1,-3){4.67}}
\put(30.09,58.33){\vector(1,-3){4.67}}
\put(29.89,58.33){\vector(1,-3){4.67}}
\put(30.19,58.33){\vector(1,-3){4.67}}
\put(29.79,58.33){\vector(1,-3){4.67}}
\put(32.49,53.33){\makebox(0,0)[lc]{$\wt\vk_1\f^3e^{ia_2}$}}
\put(30.32,30.00){\vector(1,-3){4.67}}
\put(12.33,58.34){\vector(-1,3){3.00}}
\put(29.66,58.34){\vector(1,4){2.33}}
\put(10.33,7.00){\vector(-1,4){2.33}}
\put(32.33,7.00){\vector(1,4){2.33}}
\put(11.66,21.00){\makebox(0,0)[lc]{$\G_{0,1}$}}
\put(20.99,28.34){\makebox(0,0)[ct]{$\G_{0,2}$}}
\put(32.49,22.34){\makebox(0,0)[rc]{$\G_{0,3}$}}
\put(11.32,38.00){\makebox(0,0)[lc]{$\G_{0,5}$}}
\put(42.66,15.00){\makebox(0,0)[ct]{$\G_{0,4}$}}
\put(30.66,38.00){\makebox(0,0)[rc]{$\G_{0,6}$}}
\put(11.32,50.00){\makebox(0,0)[lc]{$\G_{1,1}$}}
\put(20.99,57.00){\makebox(0,0)[ct]{$\G_{1,2}$}}
\put(31.99,51.00){\makebox(0,0)[rc]{$\G_{1,3}$}}
\put(41.99,43.34){\makebox(0,0)[ct]{$\G_{1,4}$}}
\put(18.33,4.67){\makebox(0,0)[ct]{$a)$}}
\put(62.33,17.00){\vector(1,3){4.33}}
\put(62.43,17.00){\vector(1,3){4.33}}
\put(62.23,17.00){\vector(1,3){4.33}}
\put(62.53,17.00){\vector(1,3){4.33}}
\put(62.13,17.00){\vector(1,3){4.33}}
\put(64.33,22.00){\makebox(0,0)[rc]{$s^k\wt\vk_1\f^2e^{ia_2}$}}
\put(66.67,30.33){\vector(1,0){17.33}}
\put(66.67,30.23){\vector(1,0){17.33}}
\put(66.67,30.43){\vector(1,0){17.33}}
\put(66.67,30.13){\vector(1,0){17.33}}
\put(66.67,30.53){\vector(1,0){17.33}}
\put(67.43,30.43){\makebox(0,0)[lb]{$s^k\!\f\!(\f^2\!-\!1)e^{ia}$}}
\put(88.67,16.67){\vector(1,0){16.67}}
\put(88.67,16.77){\vector(1,0){16.67}}
\put(88.67,16.57){\vector(1,0){16.67}}
\put(88.67,16.87){\vector(1,0){16.67}}
\put(88.67,16.47){\vector(1,0){16.67}}
\put(99.33,17.00){\makebox(0,0)[cb]{$\wt\vk_1\f$}}
\put(83.63,30.00){\vector(1,3){5.00}}
\put(83.73,30.00){\vector(1,3){5.00}}
\put(83.53,30.00){\vector(1,3){5.00}}
\put(83.83,30.00){\vector(1,3){5.00}}
\put(83.43,30.00){\vector(1,3){5.00}}
\put(86.33,38.00){\makebox(0,0)[lc]{$-\wt\vk_2e^{-ia_1}$}}
\put(67.00,30.00){\vector(-1,3){5.00}}
\put(67.10,30.00){\vector(-1,3){5.00}}
\put(66.90,30.00){\vector(-1,3){5.00}}
\put(67.20,30.00){\vector(-1,3){5.00}}
\put(66.80,30.00){\vector(-1,3){5.00}}
\put(65.33,35.00){\makebox(0,0)[rc]{$s^k\wt\vk_2\f e^{ia}$}}
\put(62.00,45.33){\vector(1,3){4.33}}
\put(66.33,58.67){\vector(1,0){17.33}}
\put(88.33,45.00){\vector(1,0){16.67}}
\put(88.33,44.90){\vector(1,0){16.67}}
\put(88.33,45.10){\vector(1,0){16.67}}
\put(88.33,44.80){\vector(1,0){16.67}}
\put(88.33,45.20){\vector(1,0){16.67}}
\put(97.33,45.67){\makebox(0,0)[cb]{$-\wt\vk_2\f$}}
\put(83.66,58.67){\vector(1,-3){4.67}}
\put(84.33,30.34){\vector(1,-3){4.67}}
\put(84.43,30.34){\vector(1,-3){4.67}}
\put(84.23,30.34){\vector(1,-3){4.67}}
\put(84.53,30.34){\vector(1,-3){4.67}}
\put(84.13,30.34){\vector(1,-3){4.67}}
\put(87.33,24.00){\makebox(0,0)[lc]{$\wt\vk_1e^{-ia_1}$}}
\put(65.67,21.33){\makebox(0,0)[lc]{$\G_{0,1}$}}
\put(75.00,28.67){\makebox(0,0)[ct]{$\G_{0,2}$}}
\put(86.50,22.67){\makebox(0,0)[rc]{$\G_{0,3}$}}
\put(65.33,38.33){\makebox(0,0)[lc]{$\G_{0,5}$}}
\put(96.67,15.33){\makebox(0,0)[ct]{$\G_{0,4}$}}
\put(84.67,38.33){\makebox(0,0)[rc]{$\G_{0,6}$}}
\put(65.33,50.33){\makebox(0,0)[lc]{$\G_{1,1}$}}
\put(75.00,57.33){\makebox(0,0)[ct]{$\G_{1,2}$}}
\put(86.00,51.33){\makebox(0,0)[rc]{$\G_{1,3}$}}
\put(96.00,43.67){\makebox(0,0)[ct]{$\G_{1,4}$}}
\put(66.33,58.67){\vector(-1,4){2.33}}
\put(83.67,58.67){\vector(1,4){2.33}}
\put(64.67,6.67){\vector(-1,4){2.67}}
\put(86.67,6.34){\vector(1,4){2.67}}
\put(73.00,4.67){\makebox(0,0)[ct]{$b)$}}
\put(116.00,6.66){\vector(1,3){4.33}}
\put(116.10,6.66){\vector(1,3){4.33}}
\put(115.90,6.66){\vector(1,3){4.33}}
\put(116.20,6.66){\vector(1,3){4.33}}
\put(115.80,6.66){\vector(1,3){4.33}}
\put(117.00,11.66){\makebox(0,0)[rc]{$1$}}
\put(120.34,19.99){\vector(1,0){17.33}}
\put(142.34,6.00){\vector(1,0){16.67}}
\put(142.34,6.10){\vector(1,0){16.67}}
\put(142.34,5.90){\vector(1,0){16.67}}
\put(142.34,6.20){\vector(1,0){16.67}}
\put(142.34,5.80){\vector(1,0){16.67}}
\put(153.00,6.66){\makebox(0,0)[cb]{$\f e^{i(a+a_1)}$}}
\put(137.30,19.66){\vector(1,3){5.00}}
\put(137.40,19.66){\vector(1,3){5.00}}
\put(137.20,19.66){\vector(1,3){5.00}}
\put(137.50,19.66){\vector(1,3){5.00}}
\put(137.10,19.66){\vector(1,3){5.00}}
\put(141.34,27.99){\makebox(0,0)[lc]{$-e^{ia}$}}
\put(120.67,19.66){\vector(-1,3){5.00}}
\put(120.57,19.66){\vector(-1,3){5.00}}
\put(120.77,19.66){\vector(-1,3){5.00}}
\put(120.47,19.66){\vector(-1,3){5.00}}
\put(120.87,19.66){\vector(-1,3){5.00}}
\put(117.00,28.66){\makebox(0,0)[rc]{$\f e^{ia_1}$}}
\put(138.00,20.00){\vector(1,-3){4.67}}
\put(138.10,20.00){\vector(1,-3){4.67}}
\put(137.90,20.00){\vector(1,-3){4.67}}
\put(138.20,20.00){\vector(1,-3){4.67}}
\put(137.80,20.00){\vector(1,-3){4.67}}
\put(141.00,13.66){\makebox(0,0)[lc]{$e^{ia}$}}
\put(119.34,10.99){\makebox(0,0)[lc]{$\G_{0,1}$}}
\put(128.67,18.33){\makebox(0,0)[ct]{$\G_{0,2}$}}
\put(140.17,12.33){\makebox(0,0)[rc]{$\G_{0,3}$}}
\put(119.00,27.99){\makebox(0,0)[lc]{$\G_{0,5}$}}
\put(150.34,4.99){\makebox(0,0)[ct]{$\G_{0,4}$}}
\put(138.34,27.99){\makebox(0,0)[rc]{$\G_{0,6}$}}
\put(149.67,33.33){\makebox(0,0)[ct]{$\G_{1,4}$}}
\put(142.33,35.00){\vector(1,0){15.00}}
\put(142.33,35.10){\vector(1,0){15.00}}
\put(142.33,34.90){\vector(1,0){15.00}}
\put(142.33,35.20){\vector(1,0){15.00}}
\put(142.33,34.80){\vector(1,0){15.00}}
\put(150.33,36.00){\makebox(0,0)[cb]{$-\f e^{i(a+a_1)}$}}
\put(115.67,34.66){\vector(1,3){2.00}}
\put(126.67,0.00){\makebox(0,0)[ct]{$c)$}}
\put(116.34,48.99){\vector(1,0){17.33}}
\put(116.34,48.89){\vector(1,0){17.33}}
\put(116.34,49.09){\vector(1,0){17.33}}
\put(116.34,48.79){\vector(1,0){17.33}}
\put(116.34,49.19){\vector(1,0){17.33}}
\put(125.10,49.59){\makebox(0,0)[cb]{$1$}}
\put(133.30,48.66){\vector(1,3){5.00}}
\put(133.40,48.66){\vector(1,3){5.00}}
\put(133.20,48.66){\vector(1,3){5.00}}
\put(133.50,48.66){\vector(1,3){5.00}}
\put(133.10,48.66){\vector(1,3){5.00}}
\put(137.00,56.66){\makebox(0,0)[lc]{$\f e^{ia_2}$}}
\put(116.67,48.66){\vector(-1,3){5.00}}
\put(116.77,48.66){\vector(-1,3){5.00}}
\put(116.57,48.66){\vector(-1,3){5.00}}
\put(116.87,48.66){\vector(-1,3){5.00}}
\put(116.47,48.66){\vector(-1,3){5.00}}
\put(114.00,53.66){\makebox(0,0)[rc]{$-1$}}
\put(138.00,63.66){\vector(1,0){16.67}}
\put(138.00,63.56){\vector(1,0){16.67}}
\put(138.00,63.76){\vector(1,0){16.67}}
\put(138.00,63.46){\vector(1,0){16.67}}
\put(138.00,63.86){\vector(1,0){16.67}}
\put(147.00,64.33){\makebox(0,0)[cb]{$e^{ia}$}}
\put(124.67,47.33){\makebox(0,0)[ct]{$\G_{0,2}$}}
\put(115.00,56.99){\makebox(0,0)[lc]{$\G_{0,5}$}}
\put(134.34,56.99){\makebox(0,0)[rc]{$\G_{0,6}$}}
\put(115.00,65.99){\makebox(0,0)[lc]{$\G_{1,1}$}}
\put(135.67,65.99){\makebox(0,0)[rc]{$\G_{1,3}$}}
\put(145.67,62.33){\makebox(0,0)[ct]{$\G_{1,4}$}}
\put(115.33,43.33){\vector(1,4){1.33}}
\put(133.67,49.00){\vector(1,-3){1.67}}
\put(111.67,64.00){\vector(1,3){1.33}}
\put(136.67,68.67){\vector(1,-3){1.67}}
\put(124.67,43.33){\makebox(0,0)[ct]{$d)$}}
\put(117.00,1.00){\vector(-1,4){1.33}}
\put(141.33,1.00){\vector(1,4){1.33}}
\put(140.33,39.67){\vector(1,-3){1.67}}
\end{picture}
\caption{The supports of eigenfunctions:
 a) $\p^{(0,1)}$, b) $\p^{(0,2)}$ for $(s_k,\f^2)\ne (0,1)$,
 c) $\p^{(0,1)}$, d) $\p^{(0,2)}$ for $(s_k,\f^2)=(0,1)$;
 $\p_{n,j}^{(0,\n)}=C_{n,j}^{(0,\n)}\vp$, where each $C_{n,j}^{(0,\n)}$
 is written near the  corresponding edge}
\lb{sef}
\end{figure}

\no {\bf Proof.}
Proof of i) repeats the
arguments from \cite{KL}.

\no ii) We will define the eigenfunctions of $H_k^a$ (see Fig.\ref{sef}).
If $(s_k,\f^2)\ne(0,1)$, then $\wt\vk_1\ne 0$, $\wt\vk_2\ne 0$ and the function $\p^{(0,1)}$ is given by
\[
\lb{T2-1}
\p^{(0,1)}_{n,j}=0 \ \text{for\ all}\ (n,j)\in (\Z\sm \{0,1\})\ts\N_6,\ \
\text{and} \
\p^{(0,1)}_{0,j}=0,\ j=1,3,4,\ \
\p^{(0,1)}_{1,j}=0,\ j=5,6,
$$
$$
\p^{(0,1)}_{1,4}=\f^2(\f^2-1)e^{ia}\vp,\ \
{\p^{(0,1)}_{1,1}\/\f e^{-ia_1}}={\p^{(0,1)}_{1,2}\/\f^2}={\p^{(0,1)}_{1,3}\/\f^3e^{ia_2}}
=\wt\vk_1\vp,\ \
\p^{(0,1)}_{0,5}=-\p^{(0,1)}_{0,2}=-{\p^{(0,1)}_{0,6}\/\f e^{ia_2}}=\wt\vk_2\vp,
\]
the function $\p^{(0,2)}$ is given by
\[
\lb{T2-7}
\p^{(0,2)}_{n,j}=0, \ \text{for\ all}\ (n,j)\in (\Z\sm \{0,1\})\ts\N_6,\ \
\text{and} \ \
\p^{(0,2)}_{1,j}=0,\ j\ne 4,\ \ \
\p^{(0,2)}_{0,2}=s^k\f(\f^2-1)e^{ia}\vp,
$$
$$
{\p^{(0,2)}_{0,4}\/\f}={\p^{(0,2)}_{0,3}\/e^{-ia_1}}={\p^{(0,2)}_{0,1}\/s^k\f^2e^{ia_2}}
=\wt\vk_1\vp,\qq
-{\p^{(0,2)}_{1,4}\/\f}=-{\p^{(0,2)}_{0,6}\/e^{-ia_1}}={\p^{(0,2)}_{0,5}\/s^k\f e^{ia}}
=\wt\vk_2\vp.
\]
If $(s_k,\f^2)=(0,1)$, then $\wt\vk_1=\wt\vk_2=0$ and the functions $\p^{(0,\n)}$ are given by
\[
\lb{T2-2}
\p^{(0,\n)}_{n,j}=0,\  (n,j)\in (\Z\sm\{ 0,1\})\ts\N_6,\qq
\p^{(0,\n)}_{1,j}=0, \ j\in \N_6\sm\{4\},
$$
$$
\p^{(0,1)}_{0,1}=\f e^{-ia_1}\p^{(0,1)}_{0,5}=-e^{-ia}\p^{(0,1)}_{0,6}
=e^{-ia}\p^{(0,1)}_{0,3}=\vp,\ \
-\p^{(0,1)}_{1,4}=\p^{(0,1)}_{0,4}=\f e^{i(a+a_1)}\vp,\ \
\p^{(0,2)}_{0,1}\!=0,
$$
$$
\p^{(0,2)}_{0,1}=\p^{(0,2)}_{0,3}=\p^{(0,2)}_{0,4}=0, \qq
\p^{(0,2)}_{0,2}=-\p^{(0,2)}_{0,5}=e^{-ia}\p^{(0,2)}_{1,4}=e^{-ia_2}\f\p^{(0,2)}_{0,6}
=\vp.
\]
Using \er{T2-1}-\er{T2-2}, we deduce that  $\p^{(0,\n)}$ satisfy
the Kirchhoff conditions \er{1K0}, \er{1K1}. Thus
$\p^{(0,\n)}$ are eigenfunctions  of $H_k^a$.
The operator $H_k^a$ is periodic, then each $\p^{(n,\n)}, (n,\n)\in \Z\ts
\N_2$ is an eigenfunction.
We will show that the sequence $\p^{(n,\n)}, (n,\n)\in \Z\ts
\N_2$ forms a basis for $\cH_k^a(\l)$.
Let
$h=\sum_{(n,\n)\in\Z\ts\N_2}\a_{n,\n}\p^{(n,\n)}=0.$
Then \er{T2-1}, \er{T2-7} imply for $(s_k,\f^2)\ne(0,1)$
$$
0=h|_{\G_{n,5}}=(\a_{n,1}+\a_{n,2}s^k\f e^{ia})\wt\vk_2\vp,\qq
0=h|_{\G_{n,6}}=-(\a_{n,1}\f e^{ia}+\a_{n,2})\wt\vk_2e^{-ia_1}\vp,\qq n\in\Z,
$$
which yields $\a_{n,\n}=0$ for all $(n,\n)\in\Z\ts\N_2$.
Hence $\p^{(n,\n)}$ are linearly independent.
The similar arguments show that $\p^{(n,\n)}$ are linearly independent
for $(s_k,\f^2)=(0,1)$.

For any $f\in \cH_k^a(\l)$ we will show the identity
\er{T2-3}, i.e.,
\[
\lb{eit2-2}
f=\wh f, \ \  \text{where} \  \wh f=\!\!\!\!\sum_{(n,\n)\in \Z\ts
\N_2}\!\!\!\!\wh f_{n,\n}\p^{(n,\n)}\qq
\text{and}\qq \wh f_{n,\n}\qq \text{are\ given\ by\ \er{T2-4}}.
\]
>From $\l\in\s_D$, we deduce that $f|_{\G_{n,j}}=f_{n,j}'(0)\vp$.
If $(s_k,\f^2)\ne(0,1)$, then identities  \er{T2-4}-\er{T2-7} provide for all $n\in\Z$
\begin{multline}
\lb{idFf}
\wh f|_{\G_{n,5}}
=\sum_{(\ell,\n)\in\Z\ts\N_2}\wh f_{\ell,\n}\p^{(\ell,\n)}|_{\G_{n,5}}
=(\wh f_{n,1}+\wh f_{n,2}s^k\f e^{ia})\wt\vk_2\vp=f_{n,5}'(0)\vp=f|_{\G_{n,5}},
\\
\wh f|_{\G_{n,6}}
=\!\!\!\!\sum_{(\ell,\n)\in\Z\ts\N_2}\!\!\!\!\wh f_{\ell,\n}\p^{(\ell,\n)}|_{\G_{n,6}}
=-(\wh f_{n,1}\f e^{ia}+\wh f_{n,2} )\wt\vk_2e^{-ia_1}\vp=f_{n,6}'(0)\vp=f|_{\G_{n,6}}.
\end{multline}
If $(s_k,\f^2)=(0,1)$, then identities  \er{fna},\er{T2-2} give for all $n\in\Z$
\[
\lb{T2d}
\wh f|_{\G_{n,1}}=\wh f_{n,1}\vp=f_{n,1}'(0)\vp=f|_{\G_{n,1}},\qq
\wh f|_{\G_{n,2}}=\wh f_{n,2}\vp=f_{n,2}'(0)\vp=f|_{\G_{n,2}}.
\]
Identities \er{idFf}, \er{T2d} yield
$\sum |\wh f_{n,\n}|^2<\iy $ and $\wh f\in L^2(\G^{1})$, since $f\in L^2(\G^{1})$.

Note that $\wh f$ satisfies the Kirchhoff conditions \er{1K0}, \er{1K1}
and $-\wh f_\a''+q\wh f_\a=\l \wh f_\a, \a\in \cZ_1$.
Consider the function $u=f-\wh f$.
The function $u=0$ at all vertices of $\G^1$ and then
$u_{n,j}=C_{n,j}\vp, (n,j)\in\Z\ts\N_6$.
If $(s_k,\f^2)\ne(0,1)$, then \er{idFf} gives $C_{n,5}=C_{n,6}=0$.
Let $n\in\Z$ and assume that $C_{n,1}=C$. Then
the Kirchhoff boundary conditions \er{1K0}-\er{1K1} yield
$C_{n,2}=C\f e^{ia_1}$, $C_{n,3}=C\f^2e^{ia}$, $C_{n,4}=C\f^3e^{i(a+a_1)}$,
$C_{n,1}=Cs^k\f^4e^{2ia}$. Since $C_{n,1}=C$, we obtain $C=0$ and $u=0$,
which yields \er{eit2-2}.
If $(s_k,\f^2)=(0,1)$, then \er{T2d} gives $C_{n,1}=C_{n,2}=0$.
The Kirchhoff boundary conditions
\er{1K0}-\er{1K1} yield $C_{n,5}=0$, and then $C_{n,4}=0,n\in\Z$.
Assume that $C_{n,3}=C$. Then
$C_{n,6}=-C_{0,3}=-C$ and $C_{n+1,3}=-C_{n,6}=C$ for all $n\in\Z$.
Due to $u\in L^2(\G^{1})$, we have $C=0$ and $u=0$, which yields
\er{eit2-2}.

The mapping $f\to (\wh f_{n,\n})_{(n,\n)\in \Z\ts\N_2}$ is a linear
and one-to-one  mapping from $\cH_k^a$ onto $\ell^2\os \ell^2$.
Then it is a linear isomorphism. $\BBox$

\no {\bf Proof of Theorem \ref{T3}}. (i) Using the arguments
from [BBK],[KL] and identities \er{T1-3},
\er{T1-2}, \er{T1-4}
we obtain \er{DeLk}, \er{S3b}.

\no (ii), (iii) Proof of the statements
 and identities \er{T3-1} repeats the corresponding
arguments from \cite{KL}.
 Identities \er{DeLk}, \er{scs} yield
$F_{k,\n}(\cdot,-a)=F_{N-k,\n}(\cdot,a)$,
$F_{k,\n}(\cdot,a+{\pi\/N})=F_{k+1,\n}(\cdot,a)$,
$(\n,k)\in\N_2\ts\Z_N.$
Then identities \er{T3-1} yield \er{sh1}.
By Theorem \ref{T2}, each point $\l_0\in \s_{pp}(H_k^a)$ has
an infinite multiplicity.
$\BBox$

\begin{lemma} \lb{Rpa} Let $(a,k)\in\R\ts\Z_N$.
 Then there exists an integer $n_0>1$ such that

\no (i) The function
 $D_k^-=D_0^-=\det(\cM_k+I_4)$
has exactly $4n_0$ zeros, counted with multiplicities,
 in the domain $\{\l: |\sqrt{\l}|<\pi n_0\}$
and for each $n>n_0,$ exactly two zeros,  counted with
multiplicities, in each domain
$\{\l: |\sqrt\l-\pi n-{\pi\/2}\pm\arcsin{1\/3}|<{1\/3}\}$.
There are no other zeros.

\no (ii) Each function $D_k^+=\det(\cM_k-I_4),k\in\Z_N$ has
exactly $4n_0+2$ zeros,
counted with multiplicities,
 in the domain
 $\{\l: |\sqrt{\l}|<\pi n_0+{\pi\/2}\}$,
and for each $n>n_0$  exactly one simple zero
in each domain
$\{\l: |\sqrt \l-\pi n
\pm\arccos{\sqrt{5-4c_k}\/3}|<{1\/3}\}$,
exactly one simple zero
in each domain
$\{\l: |\sqrt \l-\pi n
\pm\arccos{\sqrt{5+4c_k}\/3}|<{s_k\/3}\}$, if $s_k\ne 0$,
and exactly two zeros, counted with multiplicities,
in each domain
$\{\l: |\sqrt \l-\pi n
\pm\arccos{\sqrt{5+4c_k}\/3}|<{1\/3}\}$, if $s_k= 0$.
There are no other zeros.

\no (iii) Let $c_k\ne 0,s_k\ne 0$. Then
the function $\r_k$ has exactly $2n_0$ zeros,
counted with multiplicities,
 in the domain $\{\l: |\sqrt{\l}|<\pi n_0\}$,
and for each $n>n_0$ exactly one simple real zero
in each domain
$\{\l: |\sqrt \l-(\pi n-{\pi\/2}\pm \arcsin{s_k\/3})|
<{s_k\/3}\}$.
There are no other
zeros.

\end{lemma}

\no{\bf Proof.} Proof uses the Rouch\'e theorem
and repeats the arguments from [KL].
$\BBox$

Below we need some identities for $D_k^\pm$.
Substituting \er{DeLk} into \er{S3n} we obtain
\[
\lb{Dpm}
D_k^+
=(9F^2-g_{k,1})(9F^2-g_{k,2}),\qq
D_k^-=D_0^-
=(9F^2-h_1)(9F^2-h_2)\qq\text{on}\qq\R,
\]
where
\[
\lb{nec}
g_{k,\n}=5+F_-^2+(-1)^\n 2\sqrt{F_-^2+4c_k^2},\qq
h_\n=(1+(-1)^\n|F_-|)^2,\qq (\n,k)\in\N_2\ts\Z_N.
\]

\begin{figure}
\unitlength 1.00mm \linethickness{0.1pt}
\begin{picture}(162.00,90.00)(00.00,-10.00)
\put(10.00,-5.00){\line(0,1){75.00}}
\put(05.00,00.00){\line(1,0){155.00}}
\put(161.00,-2.00){\makebox(0,0)[cc]{$\l$}}
\put(25.00,70.00){\makebox(0,0)[cc]{$9F^2(\l)$}}
\bezier{600}(15.00,70.00)(40.00,3.00)(55.00,0.00)
\bezier{600}(55.00,0.00)(67.00,2.00)(75.00,20.00)
\bezier{600}(75.00,20.00)(90.00,58.00)(95.00,55.00)
\bezier{600}(95.00,55.00)(100.00,53.00)(113.00,30.00)
\bezier{600}(113.00,30.00)(130.00,02.00)(140.00,00.00)
\bezier{600}(140.00,00.00)(150.00,02.00)(160.00,20.00)
\bezier{600}(15.00,70.20)(40.00,3.20)(55.00,0.20)
\bezier{600}(55.00,0.20)(67.00,2.20)(75.00,20.20)
\bezier{600}(75.00,20.20)(90.00,58.20)(95.00,55.20)
\bezier{600}(95.00,55.20)(100.00,53.20)(113.00,30.20)
\bezier{600}(113.00,30.20)(130.00,02.20)(140.00,00.20)
\bezier{600}(140.00,00.20)(150.00,02.20)(160.00,20.20)
\bezier{600}(15.00,70.10)(40.00,3.10)(55.00,0.10)
\bezier{600}(55.00,0.10)(67.00,2.10)(75.00,20.10)
\bezier{600}(75.00,20.10)(90.00,58.10)(95.00,55.10)
\bezier{600}(95.00,55.10)(100.00,53.10)(113.00,30.10)
\bezier{600}(113.00,30.10)(130.00,02.10)(140.00,00.10)
\bezier{600}(140.00,00.10)(150.00,02.10)(160.00,20.10)
\bezier{600}(15.00,70.30)(40.00,3.30)(55.00,0.30)
\bezier{600}(55.00,0.30)(67.00,2.30)(75.00,20.30)
\bezier{600}(75.00,20.30)(90.00,58.30)(95.00,55.30)
\bezier{600}(95.00,55.30)(100.00,53.30)(113.00,30.30)
\bezier{600}(113.00,30.30)(130.00,02.30)(140.00,00.30)
\bezier{600}(140.00,00.30)(150.00,02.30)(160.00,20.30)
\bezier{600}(5.00,60.00)(35.00,42.00)(55.00,40.00)
\bezier{600}(55.00,40.00)(65.00,40.00)(75.00,44.00)
\bezier{600}(75.00,44.00)(85.00,49.00)(95.00,50.00)
\bezier{600}(95.00,50.00)(105.00,49.00)(130.00,40.00)
\bezier{600}(130.00,40.00)(145.00,35.00)(160.00,45.00)
\put(05.00,63.00){\makebox(0,0)[cc]{$g_{k,2}$}}
\bezier{600}(5.00,60.10)(35.00,42.10)(55.00,40.10)
\bezier{600}(55.00,40.10)(65.00,40.10)(75.00,44.10)
\bezier{600}(75.00,44.10)(85.00,49.10)(95.00,50.10)
\bezier{600}(95.00,50.10)(105.00,49.10)(130.00,40.10)
\bezier{600}(130.00,40.10)(145.00,35.10)(160.00,45.10)
\bezier{600}(5.00,30.00)(25.00,20.00)(45.00,25.00)
\bezier{600}(45.00,25.00)(60.00,30.00)(75.00,27.00)
\bezier{600}(75.00,27.00)(90.00,23.00)(105.00,25.00)
\bezier{600}(105.00,25.00)(130.00,27.00)(160.00,25.00)
\put(05.00,33.00){\makebox(0,0)[cc]{$g_{k,1}$}}
\bezier{600}(5.00,30.10)(25.00,20.10)(45.00,25.10)
\bezier{600}(45.00,25.10)(60.00,30.10)(75.00,27.10)
\bezier{600}(75.00,27.10)(90.00,23.10)(105.00,25.10)
\bezier{600}(105.00,25.10)(130.00,27.10)(160.00,25.10)
\bezier{600}(5.00,50.00)(55.00,-7.00)(95.00,1.00)
\bezier{600}(95.00,1.00)(105.00,3.00)(115.00,1.00)
\bezier{600}(115.00,1.00)(140.00,-4.00)(160.00,15.00)
\put(05.00,53.00){\makebox(0,0)[cc]{$h_{2}$}}
\bezier{600}(5.00,5.00)(20.00,-8.60)(75.00,15.00)
\bezier{600}(75.00,15.00)(90.00,20.00)(105.00,14.00)
\bezier{600}(105.00,14.00)(130.00,3.00)(160.00,1.00)
\put(05.00,8.00){\makebox(0,0)[cc]{$h_{1}$}}
\bezier{600}(5.00,50.10)(55.00,-7.10)(95.00,1.10)
\bezier{600}(95.00,1.10)(105.00,3.10)(115.00,1.10)
\bezier{600}(115.00,1.10)(140.00,-4.10)(160.00,15.10)
\bezier{600}(5.00,5.10)(20.00,-8.60)(75.00,15.10)
\bezier{600}(75.00,15.10)(90.00,20.10)(105.00,14.10)
\bezier{600}(105.00,14.10)(130.00,3.10)(160.00,1.10)
\bezier{600}(40.00,-5.00)(56.00,18.50)(70.00,-5.00)
\bezier{600}(130.00,-5.00)(143.20,10.50)(155.00,-5.00)
\put(39.00,-8.00){\makebox(0,0)[cc]{$f_{k}$}}
\put(129.00,-8.00){\makebox(0,0)[cc]{$f_{k}$}}
\bezier{600}(40.00,-5.10)(56.00,18.60)(70.00,-5.10)
\bezier{600}(130.00,-5.10)(143.20,10.60)(155.00,-5.10)
\put(22.52,0.00){\line(0,1){50.00}}
\put(22.50,-3.00){\makebox(0,0)[cc]{\tiny $\l_{2,0}^{k,+}$}}
\put(35.00,0.00){\line(0,1){23.50}}
\put(33.00,-3.00){\makebox(0,0)[cc]{\tiny $\l_{1,0}^{k,+}$}}
\put(37.40,0.00){\line(0,1){19.00}}
\put(38.50,-3.00){\makebox(0,0)[cc]{\tiny $\l_{2,1}^{0,-}$}}
\put(47.50,0.00){\line(0,1){4.50}}
\put(46.00,-3.00){\makebox(0,0)[cc]{\tiny $\l_{1,1}^{0,-}$}}
\put(48.30,0.00){\line(0,1){3.70}}
\put(51.00,-3.00){\makebox(0,0)[cc]{\tiny $r_{k,1}^{-}$}}
\put(63.20,0.00){\line(0,1){3.30}}
\put(61.00,-3.00){\makebox(0,0)[cc]{\tiny $r_{k,1}^{+}$}}
\put(63.80,0.00){\line(0,1){4.10}}
\put(66.00,-3.00){\makebox(0,0)[cc]{\tiny $\l_{1,1}^{0,+}$}}
\put(71.70,0.00){\line(0,1){13.50}}
\put(73.00,-3.00){\makebox(0,0)[cc]{\tiny $\l_{2,1}^{0,+}$}}
\put(77.60,0.00){\line(0,1){26.30}}
\put(79.00,-3.00){\makebox(0,0)[cc]{\tiny $\l_{1,2}^{k,-}$}}
\put(88.15,0.00){\line(0,1){48.80}}
\put(88.00,-3.00){\makebox(0,0)[cc]{\tiny $\l_{2,2}^{k,-}$}}
\put(101.00,0.00){\line(0,1){48.90}}
\put(101.00,-3.00){\makebox(0,0)[cc]{\tiny $\l_{2,2}^{k,+}$}}
\put(115.70,0.00){\line(0,1){25.50}}
\put(115.00,-3.00){\makebox(0,0)[cc]{\tiny $\l_{1,2}^{k,+}$}}
\put(131.50,0.00){\line(0,1){5.00}}
\put(129.00,-3.00){\makebox(0,0)[cc]{\tiny $\l_{2,3}^{0,-}$}}
\put(136.45,0.00){\line(0,1){1.50}}
\put(135.50,-3.00){\makebox(0,0)[cc]{\tiny $\l_{1,3}^{0,-}$}}
\put(137.00,0.00){\line(0,1){1.00}}
\put(140.00,-3.00){\makebox(0,0)[cc]{\tiny $r_{k,2}^{-}$}}
\put(145.50,0.00){\line(0,1){2.30}}
\put(144.60,-3.00){\makebox(0,0)[cc]{\tiny $r_{k,2}^{+}$}}
\put(145.80,0.00){\line(0,1){2.50}}
\put(149.50,-3.00){\makebox(0,0)[cc]{\tiny $\l_{1,3}^{0,+}$}}
\end{picture}
\caption{Functions $9F^2,g_k^\pm,h^\pm$ and $f_k$}
\lb{lf}
\end{figure}
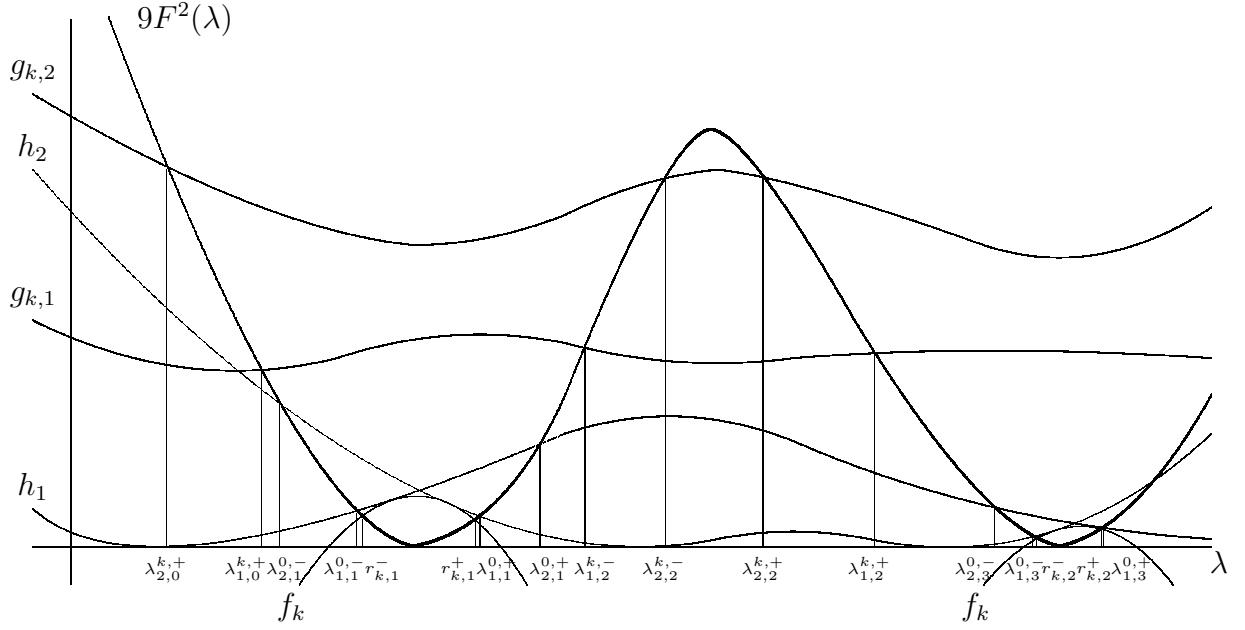

\begin{lemma}\label{T4n} Let $(a,k)\in\R\ts\Z_N$. Then

\no (i) For all $(\n,n)\in\N_2\ts\N$
the periodic and antiperiodic eigenvalues
satisfy
\[
\lb{l0}
\l_{\n,n-1}^{k,\pm}(a)=\l_{\n,n-1}^{N-k,\pm}(-a),\ \
\l_{\n,p}^{k,\pm}(a)=\l_{\n,p}^{0,\pm}(0),\ \
\g_{\n,p}^k(a)=\g_{\n,p}^0(0),\ \ \g_{1,p}^0=\vk_n,\ \ p=2n-1,
\]
\begin{multline}
\lb{b1}
\wt\l_{n-1}^+\le\l_{2,p-1}^{k,+}
\le\min\{\l_{2,p}^{0,-},\l_{1,p-1}^{k,+}\}
\le\max\{\l_{2,p}^{0,-},\l_{1,p-1}^{k,+}\}
\le\l_{1,p}^{0,-}\\ \le\e_n\le\l_{1,p}^{0,+}
\le\min\{\l_{1,p+1}^{k,-},\l_{2,p}^{0,+}\}
\le\max\{\l_{1,p+1}^{k,-},\l_{2,p}^{0,+}\}
\le\l_{2,p+1}^{k,-}\le\wt\l_{n}^-,
\end{multline}
\[
\lb{cg}
\l_{2,p}^{0,-}<\l_{1,p-1}^{k,+}\ \Leftrightarrow\
u_k(\l_{2,p}^{0,-})>0\ ;
\qqq\l_{2,p}^{0,+}>\l_{1,p+1}^{k,-}\ \Leftrightarrow\
u_k(\l_{2,p}^{0,+})>0,
\]
\[
\lb{s1a}
\cup_{n\ge 1}
\g_{\n,2n-1}^0
=\{\l\in\R:9F^2(\l)<h_\n(\l)\},\qq
\cup_{n\ge 0}\g_{\n,2n}^k
=\{\l\in\R:9F^2(\l)>g_{k,\n}(\l)\},
\]
where
$$
\g_{\n,0}^k=(-\iy,\l_{\n,0}^{k,+}),\qq
\g_{\n,n}^k=(\l_{\n,n}^{k,-},\l_{\n,n}^{k,+}),\qq
(\n,n)\in\N_2\ts\N.
$$

\no (ii) If $s_k^2(a)<s_\ell^2(a_1)$
for some $(a_1,\ell)\in\R\ts\Z_N$, then for all $n\ge 1$
\[
\lb{F2}
\l_{2,2n-2}^{k,+}(a)\!<\!\l_{2,2n-2}^{\ell,+}(a_1),\ \
\l_{1,2n-2}^{\ell,+}(a_1)\!<\!\l_{1,2n-2}^{k,+}(a),\ \
\l_{1,2n}^{k,-}(a)\!<\!\l_{1,2n}^{\ell,-}(a_1),\ \
\l_{2,2n}^{\ell,-}(a_1)\!<\!\l_{2,2n}^{k,-}(a).
\]

\end{lemma}

\no {\bf Proof.} Proof uses identities \er{Dpm} and properties
of the function $F$ and repeats the arguments from
\cite{BBK1}, Lemma 2.2 (see also Fig.\ref{lf}).
$\BBox$

\section{ Proof of Theorems \ref{Tk}-\ref{TM}}
\setcounter{equation}{0}

Let  $R_k=\{\l\in\R:\r_{k}(\l)>0\}$, if $\r_k\ne 0$, and
$R_k=\R$, if $\r_k=0$.
In Lemmas \ref{sp}-\ref{gS} we describe the set
$\s_{k,\n}=\{\l\in\R:F_{k,\n}(\l)\in[-1,1]\}$ in terms of $F$.

\begin{lemma} \lb{sp} Let $a\in\R$.
 For all $k\in\Z_N$ and $\l\in R_k$
the following relations hold true:
\[
\lb{F1+}
F_{k,\n}(\l)<1\qq \text{iff}\qq
9F^2(\l)<g_{k,\n}(\l),\qq \n=1,2,
\]
\[
\lb{g1}
F_{k,1}(\l)>-1\qq \text{iff}\qq
\{9F^2(\l)>h_1(\l)\ \text{or}\
|F_-(\l)|<c_k^2\},
\]
\[
\lb{g2}
F_{k,2}(\l)>-1\qq \text{iff}\qq
\lt\{9F^2(\l)>h_2(\l)\ \text{or}\
\{9F^2(\l) < h_1(\l)\ \text{and}\
|F_-(\l)|<c_k^2\}\rt\}.
\]
\end{lemma}

\no {\bf Proof.} Proof repeats the arguments from \cite{BBK1}, Lemma 3.1.
Using the identities \er{DeLk}, \er{nec} we write the functions
$F_{k,\n}$ in terms of $g_{k,\n},h_\n$ and prove \er{F1+}-\er{g2}.
$\BBox$

Now we describe the zeros of $\r_k$ and the functions $F_{k,\n}$
on the interval $\vk_n$.

\begin{lemma} \lb{res}
Let $c_k\ne 0$ for some $(a,k)\in\R\ts\Z_N$.
Then

\no (i) The following relation holds true:
\[
\lb{r>0}
\R\sm\ol{R_k}\ss\cup_{n\ge 1}\vk_n,\qq \text{where}\qq
\vk_n=(\l_{1,2n-1}^{0,-},\l_{1,2n-1}^{0,+}).
\]
All real zeros of $\r_k$ belong to the set $\cup_{n\ge 1}\ol\vk_n$.
Each interval $\ol\vk_n,n\ge 1$ contains
even number $\ge 0$ of zeros of $\r_k$, counted with multiplicities.

\no (ii) Let $\vk_n\not\ss R_k$. Then for each $\s=\pm$ the following
relations hold true:
\[
\lb{sFc}
\sign v_k=\const\ \ \text{on\ each}\ \
\vk_{k,n}^\s\ss R_k,
\]
\[
\lb{rg1}
\text{if}\qq v_k(\l)>0\ \text{for\ some}\ \l\in \vk_{k,n}^\s,\qq
\text{then}\qq F_{k,2}<F_{k,1}<-1\ \text{on}\   \vk_{k,n}^\s,
\]
\[
\lb{rg2}
\text{if}\qq v_k(\l)<0\ \text{for\ some}\ \l\in \vk_{k,n}^\s,\qq
\text{then}\qq -1<F_{k,2}<F_{k,1}\ \text{on}\   \vk_{k,n}^\s,
\]
where
\[
\lb{kpm}
\vk_{n,k}^-=(\l_{1,2n-1}^{0,-},r_{k,n}^-),\qq
\vk_{n,k}^+=(r_{k,n}^+,\l_{1,2n-1}^{0,+}),\qq n\ge 1
\]
If $\vk_n\ss R_k$, then
\[
\lb{un}
v_k>0\qq \text{and}\qq F_{k,2}<F_{k,1}<-1\qq\text{on}\qq  \vk_n.
\]

\no (iii) If $v_k(\l)<0$ for some $\l\in\ol{\vk_n}$, then $\r_k$
has even number $\ge 2$ of zeros on $\ol{\vk_n}$.
\end{lemma}

\no {\bf Proof.} Proof repeats the arguments from \cite{BBK1}, Lemma 3.2.
The last identity in \er{DeLk} gives
\[
\lb{rf}
\r=c_k^2(9F^2-f_k),\qq \text{where}\qq f_k=s_k^2\lt(1-{F_-^2\/c_k^2}\rt).
\]
Identities \er{nec} yield $h_1\ge f_k$.
Using the properties of $F$ we obtain \er{r>0} (see also Fig.\ref{lf}).
Using \er{rf}, Lemma \ref{sp} and the more detail analysis
(see \cite{BBK1}) we obtain \er{sFc}-\er{un} and the statement (iii).
$\BBox$

\no For each $(k,\n)\in\Z_N\ts\N_2$ we introduce the sets
\[
\lb{ss}
\gS_{k,\n}=\bigcup_{n\ge 1}\lt([\l_{\n,p-1}^{k,+},\l_{\n,p}^{0,-}]
\cup[\l_{\n,p}^{0,+},\l_{\n,p+1}^{k,-}]
\rt),\qq \gS_{k}^R=\!\!\!\!\bigcup_{\s=\pm,n\in N_\s}\!\!\!\!
\ol{\vk_{k,n}^\s},\qq N_\pm=\{n\in\N:v_k(\l_{1,p}^{0,\pm})<0\},
\]
where $p=2n-1$.
The set $\gS_{k,\n}$ is a "stable" part of
$\s_{k,\n}=\{\l\in\R:F_{k,\n}(\l)\in[-1,1]\}$, and
$\gS_{k}^R$ is an "unstable" part of $\s_{k,\n}$ (see \er{g2b}, \er{F21}).

\begin{lemma}\lb{gS} Let $(a,k,\n)\in\R\ts\Z_N\ts\N_2$.
Then the following identities hold true:
\[
\lb{cgS}
\gS_{k,\n}=\{\l\in\R:h_\n(\l)\le 9F^2(\l)\le g_{k,\n}(\l)\},
\]
\[
\lb{g2b}
\gS_{k}^R=\ca \{\l\in R_k:9F^2(\l)\le h_1(\l)\ \text{and}\
v_k(\l)\le 0\},\ \ \text{if}\ c_k\ne 0 \\
\es,\ \ \text{if}\ c_k=0\ac,
\]
\[
\lb{F21}
\s_{k,\n}=\gS_{k,\n}\cup\gS_{k}^R.
\]

\end{lemma}

\no {\bf Proof.} Proof uses the results of Lemma \ref{sp}
and repeats the arguments from \cite{BBK1}, Lemma 3.3.
Relations \er{s1a} and Lemma \ref{res} yield \er{cgS}, \er{g2b}.
Then Lemma \ref{sp} gives \er{F21}.
$\BBox$

\no{\bf Proof of Theorem \ref{Tk}.}
Substituting \er{ss} into \er{F21} and using the identity
$\s_{ac}(H_k^a)=\s_{k,1}\cup\s_{k,2}$ we obtain \er{sHk}, \er{s}.
Identities $\l_{\n,p}^{k,\pm}(a)=\l_{\n,p}^{0,\pm}(0)$
are proved in \er{l0}. Using \er{s} and the definition of $\vk_n$
in \er{kan} we obtain \er{gc}.
$\BBox$

\begin{theorem} \lb{4s}
{\rm (multiplicity of the spectrum)}
Let $(a,k,n)\in[0,{\pi\/2N}]\ts\Z_N\ts\N,p=2n-1$.
Then the following relations hold true:
\begin{multline}
\lb{eeg}
E_{2,p}^{k,\pm}=E_{2,p}^{0,\pm},\qq
E_{2,p-1}^{k,+}\le\min\{E_{2,p}^{k,-},E_{1,p-1}^{k,+}\}\le
\max\{E_{2,p}^{k,-},E_{1,p-1}^{k,+}\}
\le E_{1,p}^{k,-}\\
\le E_{1,p}^{k,+}
\le\min\{E_{1,p+1}^{k,-},E_{2,p}^{k,+}\}
\le\max\{E_{1,p+1}^{k,-},E_{2,p}^{k,+}\}\le E_{2,p+1}^{k,-},
\end{multline}
\[
\lb{mg}
E_{2,p}^{k,-}<E_{1,p-1}^{k,+}\qq \text{iff}\qq
u_k(E_{2,p}^{k,-})>0\ ;\qqq
E_{1,p+1}^{k,-}<E_{2,p}^{k,+}\qq \text{iff}\qq
u_k(E_{2,p}^{k,+})>0.
\]
Moreover,

\no (i) If $\gS_k^R\ne\es$, then
the spectrum of $H_k^a$
in $\gS_k^R$
has multiplicity 4.

\no (ii) If $\gS_k=\gS_{k,1}\cap\gS_{k,2}\ne\es$
then the spectrum of $H_k^a$ in $\gS_k$ has multiplicity 4.

\no (iii) The spectrum of $H_k^a$
in $\s_{ac}(H_k^a)\sm(\gS_k\cup\gS_k^R)$
has multiplicity 2.
\end{theorem}

\no{\bf Proof.} Relations \er{s} and \er{b1} yield \er{eeg}.
Relations \er{cg} imply \er{mg}.
Proof of the other statements repeats the arguments from \cite{BBK1},
proof of Theorem 1.2. Relations \er{rg2} provide the statement (i),
identity \er{F21} yields the statements (ii), (iii) (see \cite{BBK1}).
$\BBox$

\no {\bf Remark.}
1)
Using identities \er{s}, \er{ss} we deduce that
$\gS_k^R\ne\es$ iff
$\{E_{1,p}^{k,\s}=r_{k,n}^\s,\vk_{k,n}^\s\ne\es $
for some $(n,\s)\in\N\ts\{+,-\}\}$.

\no 2) Identities \er{s} and \er{ss} show that
$\gS_k\ne\es$ iff $E_{2,p}^{k,-}>E_{1,p-1}^{k,+}$
(or $E_{1,p+1}^{k,-}>E_{2,p}^{k,+}$) for some $n\ge 1$.
Relations \er{mg} yields $u_k(E_{2,p}^{k,-})>0$ (or
$u_k(E_{2,p}^{k,+})>0$) for this $n$.
Then the spectrum of $H_k^a$ in the interval
$(E_{1,p-1}^{k,+},E_{2,p}^{k,-})$
(or $(E_{2,p}^{k,+},E_{1,p+1}^{k,-})$)
has multiplicity 4.

Below we need the asymptotics from \cite{M}
\[
\label{Das}
F(\l)= \cos\sqrt{\l}+{O(e^{|\Im\sqrt{\l}|})\/\sqrt{|\l|}},\qqq
F_-(\l)={O(e^{|\Im\sqrt{\l}|})\/\sqrt{|\l|}}
\qqq \text{as}\qqq |\l|\to \iy.
\]

\no {\bf Proof of Theorem \ref{g}.}
(i) Estimates \er{eeg} show that the
intervals $G_{\n,k,0}=(-\iy,E_{\n,0}^{k,+})$,
$G_{\n,k,n}=(E_{\n,n}^{k,-},E_{\n,n}^{k,+})$, $(\n,n)\in\N_2\ts\N$,
satisfy:
$$
G_{1,k,0}\cap G_{2,k,m}=\es\  \text{for}\ m\not\in\{0,1\},\qq
 G_{1,k,2n-1}\cap G_{2,k,m}=\es\  \text{for}\ m\ne 2n-1,
$$
$$
 G_{1,k,2n}\cap G_{2,k,m}=\es\  \text{for}\ m\not\in\{2n-1,2n,2n+1\}.
$$
Then the gaps $G_{k,n}^a,n\ge 0$ in the spectrum $H_k^a$
satisfy
\begin{multline}
G_{k,0}^a= G_{1,k,0}\cap G_{2,k,0},\qq
G_{k,2n}^a= G_{1,k,n}\cap G_{2,k,n},\\
G_{k,4n-3}^a= G_{1,k,2n-2}\cap G_{2,k,2n-1},\qq
G_{k,4n-1}^a= G_{1,k,2n}\cap G_{2,k,2n-1},
\end{multline}
$ n\ge 1$, which yields all identities in \er{Hk1}.
Estimates \er{b1} give all inclusions in \er{Hk1}.
Relations \er{F2} give \er{gk2}.

For $(a,k)=(0,0)$ the first identity in \er{gk4}
is proved in \cite{BBK1}.
If $a\ne 0$, then $s_k\ne 0$. Let $c_k\ne 0$.
Asymptotics \er{Das} implies
$F_-(\l)\to 0$ as $\l\to +\iy$.
The last identities in \er{s} give
$E_{1,p}^{k,\pm}=r_{k,n}^\pm$, $n\ge n_0$.
By Lemma \ref{Rpa} (iii), $r_{k,n}^-<r_{k,n}^+$ for such $n$.
Hence the intervals $G_{k,4n-2}^a\ne\es$ in this case.
Let $c_k=0$. The last identities in \er{s} give
$E_{1,p}^{k,\pm}=\l_{1,p}^{0,\pm},n\ge 1$.
Identities \er{ss}, \er{cgS} yield $\l_{1,p}^{0,\pm}$ are zeros
of $9F^2-h_1$, where $h_1$ is given by \er{nec}.
We have $h_1(\l)\to 1$ as $\l\to +\iy$.
The properties of $F$ show that the function $9F^2-h_1$ has
only simple zeros on $(\l_0,+\iy)$ for sufficiently large $\l_0>0$.
Then  $\l_{1,p}^{0,-}<\l_{1,p}^{0,+}$ for all $n\ge n_0$.
Hence $G_{k,4n-2}^a\ne\es$ in this case, which yields the first
identities in \er{gk4} for $(a,k)\ne (0,0)$.

If $(a,k)\ne (0,0)$, then $s_k\ne 0$
and estimates \er{mg} give
$E_{k,4n-3}^->E_{k,4n-3}^+$ and
$E_{k,4n-1}^->E_{k,4n-1}^+$,  $n>n_0$.
Hence $G_{k,2n-1}^a=\es$ for such $k,n$, which gives
the second identities in \er{gk4}.

\no (ii) If $v_k(\l)>0$ on $\vk_n$, then
identities \er{s} show that $E_{1,p}^{k,\pm}(a)=\l_{1,p}^{0,\pm}$.
Then \er{Hk1} gives
$G_{k,4n-2}^a=(E_{1,p}^{k,-}(a),E_{1,p}^{k,+}(a))=\vk_n$.
$\BBox$

\no{\bf Proof of Theorem \ref{TM}.}
Using the results of Theorem \ref{T3} (iii) we obtain \er{TM-1}.

\no (i) Theorem \ref{g} (i) provides
$\s_{ac}(H^a)=\R\sm\cup_{n\ge 0} G_{n}^a$, where gap
$G_n^a=\cap_{k\in\Z_N}G_{k,n}^a$.
The first relations in \er{gk2} show
$G_{4n}^a=G_{0,4n}^a,n\ge 0$.
The second relations in \er{gk2} imply $G_{2n-1}^a=G_{m_1,2n-1}^a,n\ge 1$,
where $m_1={N\/2}$ or ${N-1\/2}$.
The relations $G_{k,4n-2}^a\ss\vk_n$ give $G_{4n-2}^a\ss\vk_n$.
The relations $\wt\g_{n-1}\ss G_{k,4n}^a,
\e_n\in [E_{1,2n-1}^{k,-},E_{1,2n-1}^{k,+}]$
give the corresponding relations
$\wt\g_{n}\ss G_{4n}^a,
\e_n\in [E_{1,2n-1}^-,E_{1,2n-1}^+]$.
Thus, we have proved all relations in \er{TM-3}.
Relations \er{gk2} yield \er{g3}.
Relations \er{gk4} yield \er{pgH}.

\no (ii) In order to prove asymptotics \er{aE}
we assume that $\int_0^1q(t)dt=0$.
Let $m=0$. Using the first identities in \er{s} and \er{Hk1} we obtain
$(E_{4n}^-,E_{4n}^+)=G_{4n}^a=G_{0,4n}^a=(E_{2,2n}^{0,-},E_{2,2n}^{0,+})
=(\l_{2,2n}^{0,-},\l_{2,2n}^{0,+})$.
Identities \er{nec} and Lemma \ref{Rpa} (ii)
show that $E_{4n}^\pm=\l_{2,2n}^{0,\pm}$ are zeros of the function
$9F^2-g_{0,2}$
and $\sqrt{E_{4n}^\pm}=\pi n\pm\wt\theta_0+\ve_{n}^\pm$,
where $|\ve_{n}^\pm|\le {1\/3}$ for large $n$.
Let $\l=E_{4n}^\pm$ and $\ve=\ve_{n}^\pm$.
Then \er{Das} give
$$
9F^2(\l)=9\cos^2(\pi n\pm\wt\theta_0+\ve)+O(n^{-2})
=5+4|c_0|\mp9\ve\sin(2\wt\theta_0)+O(\ve^2)+O(n^{-2}),
$$
and $F_-(\l)=O(n^{-1})$ as $n\to +\iy$. Identities \er{nec} imply
$g_{0,2}(\l)=5+4|c_0|+O(n^{-2})$.
Then $9F^2(\l)-g_{0,2}(\l)=\mp9\ve\sin(2\wt\theta_0)+O(\ve^2)+O(n^{-2})$.
Hence $\ve=O(n^{-2})$ and $\sqrt{E_{4n}^\pm}=\pi n\pm\wt\theta_0+O(n^{-2})$,
which yields  \er{aE} for $m=0$.

Let $m=1$. Identity $G_n^a=\cap_{k\in\Z_N}G_{k,n}^a$ and \er{Hk1} impliy
$E_{4n-2}^-=\max_{k\in\Z_N}E_{1,2n-1}^{k,-}$,
$E_{4n-2}^+=\min_{k\in\Z_N}E_{1,2n-1}^{k,+}$.
Identities \er{DeLk}
show that for large $n\ge 1$ the energies
$E_{1,2n-1}^{k,\pm}=r_{k,2n}^{\pm}$ are zeros of the function
$9F^2-f_k,$ where $f_k=s_k^2(1-{F_-^2\/c_k^2})$.
By Lemma \ref{Rpa} (i), (iii),
$\sqrt{E_{1,2n-1}^{k,\pm}}=\pi n-{\pi\/2}\pm\theta_k+\ve_{n}^\pm,
\theta_k=\arcsin{s_k\/3}$,
where
 $|\ve_{n}^\pm|\le {s_k\/3}$ for large $n$.
Let $\l=E_{1,2n-1}^{k,\pm}$ and $\ve=\ve_{n}^\pm$.
Then \er{Das} give
$$
9F^2(\l)=9\cos^2(\pi n-{\pi\/2}\pm\theta_k+\ve)+O(n^{-2})
=s_k^2\pm9\ve\sin(2\theta_k)+O(\ve^2)+O(n^{-2}),
$$
and $F_-(\l)=O(n^{-1})$. Then
$f_k(\l)=s_k^2+O(n^{-2})$ and
$9F^2(\l)-f_k(\l)=\pm9\ve\sin(2\theta_k)+O(\ve^2)+O(n^{-2})$.
Then $\ve=O(n^{-2})$,  which yields
$E_{1,2n-1}^{k,\pm}=(\pi n-{\pi\/2}\pm\theta_k+O(n^{-2}))^2$.
This asymptotics
and the identity $\wt\theta_1=\min_{k\in\Z_N}\theta_k$ give \er{aE} for $m=1$.

\no (iii)  Theorem \ref{g} (ii) yields the statement.
$\BBox$

\no {\bf Proof of Theorem \ref{Te}.}
(i) If $q\in L_{even}^2(0,1)$, then $F_-=0$. Let $c_k=0$.
The last identity in \er{DeLk}
implies $\r_k=0$ and identity \er{uv} yields $v_k=0$.
Using identities \er{s} $E_{1,p}^{k,\pm}=\l_{1,2n-1}^{0,\pm}$
and \er{kan}, \er{Hk1} provide $G_{k,4n-2}^a=\vk_n$.
Let $c_k\ne 0$. Then \er{uv} yields $v_k=-c_k^2<0$.
Identities \er{DeLk} show that the resonances $r_{k,n}^\pm$ are zeros of
$9F^2-s_k^2$, hence they are real.
Identities \er{s} imply $E_{1,p}^{k,\pm}=r_{k,n}^\pm$
and \er{Hk1} gives $G_{k,4n-2}^a=(r_{k,n}^-,r_{k,n}^+)$.

The identities \er{gke} for $(a,k)=(0,0)$
are proved in \cite{BBK1}.
Let $(a,k)\ne (0,0)$. Then $s_k(a)>0$ and
properties of the function $F$ yield $r_{k,n}^-< r_{k,n}^+$.
Hence $G_{k,4n-2}^a\ne\es$, which yield the first identity in \er{gke}
for this case.
Moreover, \er{uv} shows $u_k\le 0$
and estimates \er{mg} give
$E_{k,4n-3}^-\ge E_{k,4n-3}^+$ and
$E_{k,4n-1}^-\ge E_{k,4n-1}^+$,  $n\in\N$.
Hence $G_{k,2n-1}^a=\es$, which gives
the second identities in \er{gke}.
If   $s_k(a)< s_\ell(a_1)$, then the
properties of the function $F$ yield $\pm r_{\ell,n}^\mp(a_1)< \pm r_{k,n}^\mp(a)$.
Then \er{Hk1} gives $G_{k,4n-2}^a\ss G_{\ell,4n-2}^{a_1}$.

\no (ii) The first identity in \er{rge} was proved in \cite{BBKL}.
The other relations \er{rge} and \er{mfg} follow from the statement (i).
We prove \er{as1}, \er{as2}.
Recall that $E_{4n}^\pm(a)=\l_{2,2n}^{0,\pm}(a)$
are zeros of the function $9F^2-g_{0,2}$.
Since $F_-=0$, identities \er{nec} imply $g_{0,2}=5+4\cos a$.
Then $9F^2-g_{0,2}=f+4(1-\cos a),$ where $f=9F^2-9$, and
$E_{4n}^\pm(0)=\wt\l_n^\pm$ (see the first identity in \er{rge})
are zeros of the function $f$.
Let $\wt\g_n\ne\es$. Then
$$
f(\wt\l_n^\pm)=0,\qq f'(\wt\l_n^\pm)
=18F(\wt\l_n^\pm)F'(\wt\l_n^\pm)=-M_n^\pm\ne 0,\qq n\ge 1,
$$
where we have used the identities
$M_n^\pm=-F(\wt\l_n^\pm)F'(\wt\l_n^\pm)$ (see \cite{KK}).
Let $E_{4n}^\pm(a)=\wt\l_n^\pm+\ve,\ve=\ve^\pm$. Then
$f(E_{4n}^\pm(a))=\ve f'(\wt\l_n^\pm)+O(\ve^2)=-\ve M_n^\pm+O(\ve^2)$,
and using the asymptotics $4(1-\cos a)=2a^2+O(a^4),a\to 0$, we have
\[
\lb{l1}
0=9F^2(E_{4n}^\pm(a))-g_{0,2}(E_{4n}^\pm(a))=
-18 M_n^\pm\ve+2a^2+O(\ve^2)+O(a^4),\qq a\to 0,
\]
which yields $\ve=O(a^2)$. Substituting this asymptotics into \er{l1}
we obtain $\ve={a^2\/9M_n^\pm}+O(a^4)$,
which yields \er{as1}.

Let $\wt\g_n=\es$, i.e. $E\ev \wt\l_n^-=\wt\l_n^+$.
Then
$$
f(E)=f'(E)=0,\qqq f''(E)=18F(E)F''(E)=-18|F''(E)|<0.
$$
Let $E_{4n}^\pm(a)=E+\ve,\ve=\ve^\pm$.
We have
$
f(E_{4n}^\pm(a))={\ve^2\/2}f''(E)+O(\ve^3)=-9|F''(E)|\ve^2+O(\ve^3)
$
and
\[
\lb{l2}
0=9F^2(E_{4n}^\pm(a))-g_{0,2}(E_{4n}^\pm(a))=-9|F''(E)|\ve^2+2a^2+O(\ve^3)+O(a^4),
\]
which yields $\ve=O(a)$.
Substituting this asymptotics into \er{l2},
we obtain  $\ve^2={2a^2\/9|F''(E)|}+O(a^3)$,
which yields \er{as2}.

We prove \er{as3}. Using $G_{k,4n-2}^a\ss G_{\ell,4n-2}^{a_1}$
and \er{s}, \er{Hk1}  we obtain
$E_{4n-2}^\pm(a)=E_{1,2n-1}^{0,\pm}(a)=r_{0,n}^\pm(a)$.
Identities \er{DeLk} show that $r_{0,n}^\pm(a)$ are zeros
of the function $9F^2-\sin^2 a$ and $r_{0,n}^-(0)=r_{0,n}^+(0)=\e_n$.
Let $r_{0,n}^\pm(a)=\e_n+\ve,\ve=\ve_n^\pm$. Then
$0=9F^2(r_{0,n}^\pm(a))-\sin^2 a=9 (F'(\e_n))^2\ve^2- a^2+O(\ve^3)+O(a^3)$,
which yields $\ve=O(a)$.
Then $\ve=\pm {a\/3|F'(\e_n)|}+O(a^2)$, which gives \er{as3}.
$\BBox$

\begin{proposition}
\lb{pro}
Let $q=q_\ve={1\/\ve}\d(t-{1\/2}-2\ve),\ve>0$.
Then for any $n_0\in\N$ there exists $\ve_0>0$ such that
\[
\lb{del}
G_{4n-2}^a(q_\ve)=G_{4n-2}^0(q_\ve)\ne\es\qq \text{all}\qq
(a,\ve,n)\in\R\ts(0,\ve_0)\ts\N_{n_0}.
\]
\end{proposition}

\no {\bf Proof.}
We have $F_-(\l,q_\ve)
={\sin 4\ve\sqrt\l\/2\ve\sqrt\l}=2+O(\ve^2),\ve\to 0$ uniformly
on any bounded subset of $\C$
(see, for example, \cite{BBKL}).
Then $|F_-(\l,q_\ve)|\ge\cos^2a$
for each $(a,\ve,\l)\in\R\ts(0,\ve_0)\ts C_0(\wt R)$
for any $\wt R>0$ and some $\ve_0>0$,
where $C_0(\wt R)=\{\l\in\C:|\l|<\wt R\}$.
Using Theorem \ref{TM} (iii) we obtain
$G_{4n-2}^a(q_\ve)=G_{4n-2}^0(q_\ve)=\vk_n(q_\ve)$ for all
$(a,\ve,n)\in\R\ts(0,\ve_0)\ts\N_{n_0}$
for any $n_0\in\N$ and some $\ve_0>0$.
Then identities \er{gc}
yield \er{del}.
$\BBox$

In order to describe the gaps in the spectrum of
$H_k^a,H^a$ we need

\begin{Def}
Let $g=(\l_1,\l_2)$ be a gap in the spectrum of $H_k^a$ or $H^a$.

\no (i) If $\l_1,\l_2$ are zeros of $D_k^+$ $($or $D_k^- )$,
then $g$ is a periodic $($or antiperiodic$)$ gap.

\no (ii) If $\l_1,\l_2$ are zeros of $\r_k$,
then $g$ is a resonance gap.

\no (iii) If one of the numbers $\l_1,\l_2$ is a zero
of $D_k^-$
and other is a zero of $D_k^+$ $($or $\r_k )$,
then $g$ is a p-mix gap $($or r-mix gap$)$.

\end{Def}

\no In our armchair model there is no
a gap $(\l_1,\l_2)$, where one
of the numbers $\l_1,\l_2$
is a zero of $D_k^+$ and other is a zero of $\r_k$.

\begin{theorem}\lb{C} Let $(a,k,n)\in[0,{\pi\/2N}]\ts\Z_N\ts\N$. Then
(i) $G_{k,4n}^a$ are periodic gaps.

\no (ii) $G_{k,2n-1}^a$ are p-mix gaps. Let $s_k\ne 0$.
Then $G_{k,2n-1}^a=\es$, $n\ge n_0$
for some $n_0\ge 1$.

\no (iii) $G_{k,4n-2}^a$ are antiperiodic,
or resonance, or r-mix gaps.
Moreover, for some $n_0\ge 1$
$$
\ca \text{if}\ s_k=0,\ \text{then}\
G_{k,4n-2}^a\ \text{are\ antiperiodic\ gaps\ or}\ G_{k,4n-2}^a=\es,\
\text{and}\
G_{k,4n-2}^a=\es\ \text{for}\ n\ge n_0\\
\text{if}\ c_k\ne 0,s_k\ne 0,\ \text{then}\
G_{k,4n-2}^a\ \text{are resonance gaps}\
\text{for}\ n\ge n_0\\
\text{if}\ c_k= 0,\ \text{then\ all}\
G_{k,4n-2}^a\ \text{are antiperiodic gaps}.
\ac\!\!\!\!\!\!
$$
Let, in addition, $q\in L_{even}^2(0,1)$.
Then for all $n\ge 1$
$$
\ca \text{if}\ s_k=0,\ \text{then}\ E_{1,2n-1}^{k,-}=E_{1,2n-1}^{k,+}=\e_n
\ \text{and}\
G_{k,4n-2}^a=\es\\
\text{if}\ c_k\ne 0,s_k\ne 0,\ \text{then}\
G_{k,4n-2}^a\ \text{are resonance gaps}.
\ac
$$

\end{theorem}

\no {\bf Proof.} Identities \er{s}, \er{Hk1} show that
$G_{k,4n}^a$ are periodic gaps,
$G_{k,2n-1}^a$ are p-mix gaps and
$G_{k,4n-2}^a$ are antiperiodic, or resonance,
or r-mix gaps.
Asymptotics \er{Das} implies
$F_-(\l)\to 0$
as $\l\to +\iy$. Estimates \er{mg} give that
$E_{k,4n-3}^->E_{k,4n-3}^+$ and
$E_{k,4n-1}^->E_{k,4n-1}^+$ for $k\ne 0$ and large $n>1$.
Hence $G_{k,2n-1}^a=\es$ for such $k,n$.

If $s_k=0$, then \er{DeLk} gives $\r_k=9F^2 c_k^2$ and $r_{k,n}^-=r_{k,n}^+=\e_n$.
Identities \er{s} give
$E_{1,p}^{k,\pm}=\ca
\l_{1,p}^{0,\pm}\qq \text{if}\qq
v_k(\l_{1,p}^{0,\pm})\ge 0\\
\e_n\qq\ \text{if}\qq
v_k(\l_{1,p}^{0,\pm})<0
\ac\!\!\!$.
Hence
$G_{k,4n-2}^a$ are antiperiodic gaps or $G_{k,4n-2}^a=\es$.
Moreover, $v_k(\l_{1,p}^{0,\pm})<0$ and $E_{1,p}^{k,-}=E_{1,p}^{k,+}=\e_n$
for all large $n\ge 1$. Hence $G_{k,4n-2}^a=\es$ for sufficiently large $n\ge 1$.
If $c_k=0$, then $v_k(\l)\ge 0$ for all $\l\in\R$ and identities \er{s} provide
$E_{1,p}^{k,\pm}=\l_{1,p}^{k,\pm}$.
Hence $G_{k,4n-2}^a$ are antiperiodic gaps.
Let $c_k\ne 0$. Then identities \er{s} give
$G_{k,4n-2}^a$ are antiperiodic, or resonance, or r-mix gaps.
Using $F_-(\l)\to 0$ as $\l\to +\iy$ we obtain $v_k(\l)<0$ for large $\l>0$, then
$G_{k,4n-2}^a$ are resonance gaps for sufficiently large $n\ge 1$.

If $q\in L_{even}^2(0,1)$ and $c_k\ne 0$, then $F_-=0$ and
$v_k<0$.
Identities \er{s} show that
$E_{1,p}^{k,\pm}=r_{k,n}^\pm$ in this case.
The last identity in \er{DeLk} yield $\r_k=(9F^2-s_k^2)c_k^2$.
If $c_k\ne 0,s_k\ne 0$, then the
properties on the function $F$ show that $r_{k,n}^-<r_{k,n}^+$
for all $n\ge 1$.
Then  $G_{k,4n-2}^a=(r_{k,n}^-,r_{k,n}^+),n\ge 1$ are resonance gaps.
$\BBox$

\end{document}